\documentclass[10pt]{article}

\usepackage{babel}
\usepackage[utf8]{inputenc}
\usepackage[shortlabels]{enumitem}
\usepackage[a4paper,bottom=4cm]{geometry}  
\usepackage{mathtools} 
\usepackage{booktabs}   
\usepackage{amsfonts}
\usepackage{amsthm}
\usepackage{graphicx}
\usepackage{subcaption}
\usepackage{amsfonts} 
\usepackage{placeins} 
\usepackage{amsmath}
\usepackage{amssymb}  
\usepackage{float}
\usepackage{url}
\usepackage{bm}
\usepackage[font=small]{caption}
\usepackage[]{fancyhdr}  
\usepackage{mathtools}
\usepackage{booktabs} 
\usepackage{soul}
\usepackage[dvipsnames,table]{xcolor}
\usepackage[normalem]{ulem}

\usepackage[colorlinks=true,breaklinks=true,bookmarks=true,urlcolor=myblue, citecolor=myblue,linkcolor=myblue,bookmarksopen=false,draft=false]{hyperref}

\definecolor{myblue}{RGB}{0, 0, 160}

\theoremstyle{plain}
\newtheorem{theorem}{Theorem}[section]
\newtheorem{lemma}[theorem]{Lemma}

\newtheorem{proposition}[theorem]{Proposition}

\newtheorem{assumption}[theorem]{Assumption}
\theoremstyle{definition}
\newtheorem{definition}[theorem]{Definition}

\theoremstyle{remark}
\newtheorem{remark}[theorem]{Remark}

\title{Parametric Semidefinite Programming:\\ Geometry of the Trajectory of Solutions}

\author{Antonio Bellon\thanks{Faculty of Electrical Engineering, Czech Technical University in Prague, Karlovo náměstí 13, Prague 121 35, Czech Republic}
\and Vyacheslav Kungurtsev${}^*$
\and Jakub Mareček${}^*$ 
\and Didier Henrion${}^*$\thanks{LAAS-CNRS,
7 avenue du Colonel Roche,
31 400 Toulouse,
France} 
}
\date{}

\begin{document}
\newgeometry{bottom=4.2cm}
\maketitle

\begin{abstract}
    In many applications, solutions of convex optimization problems are updated on-line, as functions of time. In this paper, we consider \textit{{parametric}} semidefinite programs, which are linear optimization problems in the semidefinite cone whose coefficients (input data) depend on \textit{{a time parameter}}. We are interested in the geometry of the solution (output data) trajectory, defined as the set of solutions depending on the \textit{{parameter}}.
    We propose an exhaustive description of the geometry of the solution trajectory.
    As our main result, we show that  only six distinct behaviors can be observed at a neighborhood of a given point along the solution trajectory. Each possible behavior is then illustrated by an example. 
\end{abstract}
  
\section{Introduction}

A \textit{semidefinite program} (SDP) is a convex constrained optimization problem wherein one wants to optimize a linear objective function over the intersection of the cone of positive semidefinite matrices with an affine space. In this paper we consider \textit{{parametric} SDPs}, which are problems of the form 
\[
\label{primal_sdp}
\tag{P$_t$}
\begin{aligned}
\min_{X\in\mathbb{S}^n}\ &\langle C(t),X\rangle\\
s.t. \ \ & \mathcal{A}(t)[X]=b(t)\\
&X\succeq0
\end{aligned}
\]
whose coefficients depend on a parameter $t$ belonging to a given open interval $T=(t_{i}, t_{f}) \subseteq\mathbb{R}$.
\newline

The goal of \eqref{primal_sdp} is to optimize a linear objective function over a feasible region defined by non-linear constraints, {where the data of the problem depend on the parameter $t$, which we will often refer to as time}. The objective is to minimize the scalar product $\langle C(t),X\rangle$ between two matrices of $\mathbb{S}^n$, the vector space of symmetric matrices of size $n$ with real entries. 
The {parametric}  feasible region is an intersection of the semidefinite cone $\mathbb{S}^n_+=\{X\in\mathbb{S}^n\ |\ v^TXv\ge0,\ \forall v\in \mathbb{R}^n\}$ with a {parametric} affine subspace described by linear equations.  
The notation $X\succeq0$ is a shortcut for $X \in \mathbb{S}^n_+$. The notation $\mathcal{A}(t)[X]=b(t)$ models linear equations that $X$ must also satisfy: $\langle A_i(t),X\rangle=b_i(t)$ for $i=1,\dots,m,$ where $A_i(t) \in \mathbb{S}^n$ are given matrices and $b_i(t)$ are given scalars { depending on the parameter $t$}. 
Thus, problem $(P_t)$ is convex and its feasible region at any value of parameter $t\in T$ is an affine section of the semidefinite cone, often referred to as a \textit{spectrahedron}. 
In Section \ref{sec: prelim}, we present our notation in more detail.
\newline 
\restoregeometry
Parametric SDP appears in numerous applications, where the parameter often coincides with time.
For example, in power systems, semidefinite programming relaxations of the so-called alternating-current optimal power flow (ACOPF) are very successful, cf.  \cite{lavaei2011zero}.
Tracking of a trajectory of solutions to ACOPF with supply and demand varying over time is crucial for a transmission system operator, who decides on the activation of ancillary services to balance the transmission system, cf. \cite{liu2018coordinate}. 
In general, our goal is to understand certain properties of such solution trajectories, which would make it possible to design algorithms for {parametric} SDP with guarantees on their performance.

\section*{Background and Contribution}

SDP can be thought of as a generalization of linear programming (LP) with a number of applications in data science. 
\cite{anjos2011handbook} offered a snapshot of the state of the art in the areas of SDP, conic optimization, and polynomial optimization. 
Polynomial optimization problems can be approximated via a hierarchy of SDP problems of increasing size, developed by  \cite{lasserre2001global}, also known as the moments - sum of squares (SOS) hierarchy. Many problems in control theory can be reduced to solving polynomial equations,
polynomial inequalities, or polynomial differential equations, and they can hence be often solved approximately by the moment-SOS hierarchy, see \cite{hkl20} for a recent overview. Applications in theoretical computer science include  approximation algorithms for fundamental problems like the Max-Cut problem or coloring problems, quantum information theory, robust learning, and estimation.
\newline
\par
The geometry of SDP, that is, the geometry of the feasible region of an SDP problem, is well understood. We refer to \cite[Chapter 3 by Pataki]{wolkowicz2012handbook} for an excellent overview. Likewise,  solution regularity (duality, strict feasibility, uniqueness of the solution, strict complementarity, non-degeneracy) and its prerequisites are well understood; see for example \cite{alizadeh1997complementarity} where the relation between uniqueness of the solution, non-degeneracy of the solutions and strict complementarity is discussed. 
\newline

\par 
Here, our purpose is to study the behavior of the trajectory of the solutions to {parametric} SDP. Around points of the trajectory satisfying strict complementarity and uniqueness, by means of the implicit function theorem, one can show that the trajectory defines a smooth curve (Theorem \ref{thm: single_valued_differentiable}). When this fails to happen, a number of irregular behaviors may arise. The main result of this paper (Theorem \ref{thm: main_result}) consists of a complete classification of such points. So far, to the best of our knowledge, a complete classification of types of  behavior of points making up the trajectory of solutions has not been proposed.
Here, we suggest one based on a purely logical construction, whose definitions use set-valued analysis.
In particular, we use the Painlevé-Kuratowski extension of the notion of continuity to the case of set-valued functions, so as to reason about continuity properties at values of the time parameter, when there are multiple solutions.  
Informally, we now define the types of points that our classification comprises. This is based on the geometry of the trajectory of solutions parametrized over a time interval. Before a given time, we assume that the trajectory is regular and follows a continuous curve. Then at the time of interest, we can distinguish between the following situations:
\newline

\begin{itemize}
    \item \textbf{Non-differentiable point:} the trajectory is single-valued but not differentiable;
    \item \textbf{Discontinuous isolated multiple point:} a loss of continuity causes a loss of uniqueness of the solution, implying a multiple-valued solution. After the point, uniqueness is restored, and hence the loss of uniqueness is isolated;
    \item \textbf{Discontinuous non-isolated multiple point:} a loss of continuity causes a loss of uniqueness of the solution, implying a multiple-valued solution. After the point, uniqueness is not restored hence the loss of uniqueness is not isolated; 
    \item \textbf{Continuous bifurcation point:} the trajectory splits into several distinct branches. This results in a loss of uniqueness which still preserves continuity;  
    \item \textbf{Irregular accumulation point:} accumulation point of a set made of either bifurcation points or discontinuous isolated multiple points.
    \\
\end{itemize}
The formal definitions of the point types can be found in Def. \ref{def: type_regular}, \ref{def: type_non_diff},  \ref{def: type_isol_non_uni}, \ref{def: type_non_isol_non_uni}, \ref{def: type_bif}, and \ref{def: type_accu} 
\newline

We believe that a first contribution of this paper is precisely the definition of these types of points. 
In this respect, our approach was deeply inspired by \cite[Chapter 2]{guddat1990parametric} where a classification of solutions to univariate parametric nonlinear constrained optimization problems (NLPs) is proposed. There, critical points satisfying first-order optimality conditions are considered. Under precise algebraic conditions, these points are ``non-degenerate'' (see Remark \ref{rem: non_degeneracy}). The local behavior of such points is then shown to be regular. 
If a critical point is instead ``degenerate'' then, according to which algebraic condition is satisfied, the point is classified into four different types. Our approach is the same in spirit, in that we also start by considering algebraic conditions ensuring a regular behavior. As a main difference, we classify irregular points according to the behavior of the trajectory of solutions at the point considered rather than according to different sets of algebraic conditions (see Remark \ref{rem: non_degeneracy}). 
We point out that both the following classification results consider an optimal solution $X^*$ at a time $t^*$ under the assumption that optimal solutions are unique in a sufficiently small left time neighborhood of $t^*$.
The main results that we present in this paper are Theorem \ref{thm: main_result} and Theorem \ref{thm: finite_bad_result}, which we informally state here.
 
\begin{theorem}[Informal statement of Theorem \ref{thm: main_result}]
\label{thm: informal}
Under assumptions of Linear Independence Constraint Qualification (LICQ, cf. Assum. \ref{ass: licq}), existence of strictly feasibile point (cf. Assum. \ref{ass: strict_feas}) and continuity of the data with respect to time (cf. Assum. \ref{ass: data_continuity}),  the trajectory can only be comprised of points of the six types described above. 
\end{theorem}
\begin{theorem}[Informal statement of Theorem \ref{thm: finite_bad_result}]
Under the same assumptions of Theorem \ref{thm: informal}, suppose that the problem data are polynomial functions of time and that there exists a generic non-singular time (see Def. \ref{def: non_singular_times}). Then the trajectory is comprised of only regular points (cf. Def. \ref{def: type_regular}), non-differentiable points (cf. Def. \ref{def: type_non_diff}), or isolated multiple points (cf. Def. \ref{def: type_isol_non_uni}). In other words, non-isolated discontinuous multiple points (cf. Def. \ref{def: type_non_isol_non_uni}), bifurcation (cf. Def. \ref{def: type_bif}) points, and irregular accumulation points (cf. Def. \ref{def: type_accu}) cannot appear. 
\end{theorem}
 
Notice that while we only assume continuity of the data for Theorem \ref{thm: main_result}, we need stronger regularity assumptions to guarantee the validity of Theorem \ref{thm: finite_bad_result} (as well as Theorems \ref{thm: single_valued_differentiable} and \ref{thm: finiteness_of_singularity}).
\newline
\par
Other than interesting for a purely theoretical study, we believe that these results could be useful for algorithmic design as follows. If one can guarantee that the conditions of Theorem \ref{thm: finite_bad_result} are satisfied, algorithms for {parametric} optimization need not consider the behaviors corresponding to Definitions \ref{def: type_bif}, \ref{def: type_non_isol_non_uni}, and \ref{def: type_accu}. 
If, however, one would like to develop a solver for the case where only Assumptions \ref{ass: licq}, \ref{ass: strict_feas}, and \ref{ass: data_continuity} are satisfied, some rather pathological behaviors, such as  non-isolated discontinuous multiple points (Def. \ref{def: type_non_isol_non_uni}) or bifurcation points (Def. \ref{def: type_bif}), need to be to considered. In this respect, we believe that our work has the merit of clarifying and making explicit the nature of the irregularities of the trajectories to {parameteric} SDP.
Even though the precise algorithmic consequences will clearly be strongly dependent on the type and the properties of the algorithm in use, hence lying beyond the boundaries of our discussion, we very much hope that our study leads to the development and the improvement of practical algorithms for {parameteric} semidefinite programming.

\section*{Previous Work}
 
Following the pioneering contribution of \cite{goldfarb1999parametric}, who first studied the properties of the optimum as a function of a varying
parameter and extended the concept of the optimal partition from LP to SDP, a number of important papers appeared recently.
In decreasing order of generality of the dependence of problem data (coefficients) on the parameter (time):
\vspace{0.3cm}

$\bullet$ Al-Salih and Bohner \cite{al2018linear} studied LP on time scales, which allows for the mixing of difference and differential operators in a broad class of extensions of LP models. While very elegant, mathematically, it seems non-trivial to extend this approach to SDP;

\vspace{0.2cm}
$\bullet$ Wang, Zhang, and Yao \cite{wang2009separated} studied a broad family of parametric optimization problems, which are known as separated continuous conic programming (SCCP). They  developed a
strong duality theory for SCCP and proposed a polynomial-time approximation algorithm that solves an SCCP to any
required accuracy.
This algorithm does not, however, seem easy to extend to SDP; 
\vspace{0.2cm}

$\bullet$ Mohammad-Nezhad \cite{mohammad2019conic}, together with Terlaky 
\cite{mohammad2020parametric}, Haunstein, and Tang \cite{hauenstein2019computing} are perhaps the closest to our work, in spirit. 
Their dependence of the problem on the data is assumed to be linear, which is a more restrictive assumption than the one we use. Moreover, they do not provide a complete characterization of the possible behaviors of the trajectory of the solutions.
In part, we build upon their theoretical results, but instead of building upon the concepts of non-linearity intervals, invariancy intervals, and transition points, we use a purely set-valued analysis approach;

\vspace{0.2cm}
$\bullet$ El Khadir \cite{el2020semidefinite} and Ahmadi \cite{ahmadi2021time} studied {time-varying semidefinite programs (TV-SDPs)} in the setting where the data vary with known polynomials of the parameter and showed that under a strict feasibility assumption, restricting the solutions to be polynomial functions of the parameter does not change the optimal value of the TV-SDP. They also provided a sequence of SDP problems that give upper bounds on the optimal value of a TV-SDP converging to the optimal value. 
\newline

{
\par Let us remark here the difference in the literature between TV-SDPs and parametric SDPs. The~first type is the one considered by the aforementioned \cite{el2020semidefinite} and \cite{ahmadi2021time}. There, the constraints at a given time point are linked to to the solutions at the previous times via kernel terms.  In this case, the solutions are thought as measurable functions, which are required to satisfy the constraints only on a set of times that is the complement of a measure-zero set, i.e. \textit{almost everywhere}. 
Instead, we consider the easier case of parametric SDP, where constraints are independent through time and the solutions can be thought of as set-valued maps. 
This approach, considered by   \cite{goldfarb1999parametric}, \cite{yildirim2004} and \cite{hauenstein2019computing}, simply assumes that the coefficients of the SDP depend on a parameter.
}

\section{Preliminaries}
\label{sec: prelim}
In this section, we expose the tools needed to state and prove our main result. In Subsection \ref{subsec: properties_of_SDP} we first review geometric properties of SDP. In Subsection \ref{subsec: set_valued_analysis}  we survey continuity properties of the optimal and feasible sets of {parametric} SDP, considered as set-valued maps, in terms of inner and outer semi-continuity and Painlevé-Kuratowski continuity, to adopt the notions of \cite{rockafellar2009variational}.
Then, in Subsection \ref{subsec: Regularity_properties_of_optimal_set_map} we show that the existence of a unique pair of strictly complementary primal and dual solutions at a value of the time parameter $\hat t$ implies that there is a neighbourhood of $\hat t$ where both the primal and dual optimal trajectory have a regular behavior. Finally, we observe that under fairly weak assumptions, among which the existence of a generic non-singular point in the parameterization interval, the number of points where strict complementarity or uniqueness is lost is finite. 

\subsection{SDP optimality conditions and properties}
\label{subsec: properties_of_SDP}
 
In this section, we assume for notational simplicity that the data $\mathcal{A},b,C$ are not parameter-dependent. Let us review geometric properties of SDP in primal form
\[
\tag{P}
\begin{aligned}
\min_{X\in\mathbb{S}^n}\ &\langle C,X\rangle\\
s.t. \ \ & \mathcal{A}[X]=b \\
&X\succeq0
\end{aligned}
\]
and dual form
\[ 
\tag{D}
\begin{aligned}
\max_{y\in\mathbb{R}^m,\,Z\in\mathbb{S}^n}\ &\langle b,y\rangle\\
s.t. \quad\ \ & \mathcal{A}^*[y]+Z=C\\
&Z\succeq0.
\end{aligned}
\]
The linear operator $\mathcal{A}$ maps $X\in\mathbb{S}^n$ to $(\langle A_1,X\rangle,\dots,\langle A_m,X\rangle)\in\mathbb{R}^m$ where $A_i \in \mathbb{S}^n$ are given matrices for $i=1,\dots,m$ and $b\in\mathbb{R}^m$. $\mathcal{A}^*[y]=\sum_{i=1}^m A_iy_i$ is the linear operator adjoint to $\mathcal{A}$.
\newline
\begin{remark}
\label{rem: LICQ}
Throughout this paper, we assume that the linear operator $\mathcal{A}$ is surjective (see Assum. \ref{ass: licq}). Then, given a matrix $Z\in \mathbb{S}^n$ satisfying the dual constraint $\mathcal{A}^*[y]+Z=C$ for some $y\in\mathbb{R}^m$, $y$ can be uniquely determined by solving the linear system $(\mathcal{A}\mathcal{A}^*)[y]=\mathcal{A}(C-Z)$. We exploit this fact and when discussing a dual point $(y,Z)$ we will often omit $y$ and refer to a dual point simply as a matrix $Z\in\mathbb{S}^n$. 
\\
\end{remark}

\par
We call a matrix $X$ satisfying the constraints of $(P)$ a \textit{primal feasible point}, a matrix $Z$ satisfying the constraints of $(D)$ a \textit{dual feasible point}, a pair of matrices $(X,Z)$ satisfying the constraints of $(P,D)$ a \textit{primal-dual feasible point}.
We call a solution $X^*$ to $(P)$ a \textit{primal optimal point}, a solution $Z^*$ to $(D)$ a \textit{dual optimal point}, a solution $(X^*,Z^*)$ to $(P,D)$ a \textit{primal-dual optimal point}.
First, we recall that for the primal-dual pair of SDPs $(P,D)$ a set of first order optimality sufficient conditions is available.
Given two matrices $X$ and $Z$ of $\mathbb{S}^n$, their scalar product is denoted by $\langle X, Z\rangle=\operatorname{trace}(XZ)=\sum_{i,j=1}^nX_{i,j}Z_{i,j}$.
\newline

\par
\begin{definition}[KKT conditions]
A primal-dual feasible point $(X,Z) \in\mathbb{S}^n\times\mathbb{S}^n$ satifies the \textit{Karush-Kuhn-Tucker conditions} (KKT) for $(P,D)$ if
\begin{equation}
\label{eq: KKT}
\begin{array}{l}
\mathcal{A}[X]=b\\ \mathcal{A}^*[y]+Z=C \\
X, Z \succeq 0 \\
\langle X,Z\rangle=0
\end{array}
\tag{KKT}
\end{equation}
for some $y\in\mathbb{R}^m$. 
\\
\end{definition}

\par
It is well-known that for a convex optimization problem, conditions \eqref{eq: KKT} are sufficient for optimality. Under strict feasibility, the KKT conditions are also necessary.
\\

\begin{definition}[Strict feasibility]
We say that \textit{strict feasibility} holds for $(P)$ (or that $(P)$ is \textit{strictly feasible}) if there exists an interior point of the primal feasible region. That is, there exists a matrix $X\succ0$ satisfying $\mathcal{A}[X] = b$. Similarly, strict feasibility holds for $(D)$ (or $(D)$ is \textit{strictly feasible}) if there exists an interior point of the dual feasible region. That is, there exist $y\in\mathbb{R}^m$ and a matrix $Z\succ 0$ satisfying $\mathcal{A}^*[y]+Z=C$.
\\
\end{definition}

\begin{definition}[Strict complementarity]
A primal-dual optimal point $(X,Z)$ is said to be \textit{strictly complementary} if $\operatorname{rank}(X)+\operatorname{rank}(Z)=r+s=n$. A primal-dual problem $(P,D)$ satisfies \textit{strict complementarity} if there exists a strictly complementary primal-dual optimal point $(X,Z)$.
\\
\end{definition}

We now introduce the definitions of primal and dual non-degeneracy. All the definitions and results exposed below until the end of the subsection are due to \cite{alizadeh1997complementarity}. 
\newline

\par
\begin{definition}[Primal non-degeneracy]
\label{def: primal_non_deg}
We say that a primal feasible point $X$ is \textit{primal non-degenerate} if
\begin{equation*}
    \mathcal{N}(\mathcal{A})+\mathcal{T}_{X}=\mathbb{S}^{n},
\end{equation*}
where $\mathcal{N}(\mathcal{A})=\left\{Y \in \mathbb{S}^{n} | \langle A_{i}, Y\rangle =0 \text { for all } i=1,\dots,m\right\}$,
\[
    \mathcal{T}_{X}=\left\{Q\left(\begin{array}{cc}{U} & {V} \\ {V^{T}} & {0}\end{array}\right) Q^{T}\ \Big |\ U \in \mathbb{S}^{r}, V \in \mathbb{R}^{r \times (n-r)}\right\}
\]
is the tangent space at $X$ in $\mathbb{S}^n_+$ with $r=\operatorname{rank}(X)$, $Q=Q^T\in \mathbb{R}^{n\times n}$ is an orthogonal matrix such that its columns form a basis of eigenvectors relative to the eigenvalues $\lambda_i$ of $X$:
\begin{equation}
\label{eq: primal_diagonal}
X = Q\operatorname{diag}(\lambda_{1},\dots,\lambda_{r},0,\dots,0)Q^{T}.
\end{equation}

\end{definition} 
\begin{definition}[Dual non-degeneracy]
\label{def: dual_non_deg}
We say that a dual feasible point $Z$ is \textit{dual non-degenerate} if
    \begin{equation*}
    \mathcal{R}(\mathcal{A})+\mathcal{T}_{Z}=\mathbb{S}^{n},
    \end{equation*}
where $\mathcal{R}(\mathcal{A})=\operatorname{span}( A_{1},\dots,A_{m})$ and
\begin{equation*}
    \mathcal{T}_{Z}=\left\{\widetilde{Q}\left(\begin{array}{cc}{0} & {V} \\ {V^{\mathrm{T}}} & {W}\end{array}\right) \widetilde{Q}^{\mathrm{T}}\ \Big |\ W \in \mathbb{S}^{\mathrm{s}},\ V \in \mathbb{R}^{(n-s) \times s}\right\}
\end{equation*}
is the tangent space at $Z$ in $\mathbb{S}^n_+$, $s=\operatorname{rank}(Z)$, $\widetilde{Q}=\widetilde{Q}^T\in \mathbb{R}^{n\times n}$ is an orthogonal matrix such that its columns form a basis of eigenvectors relative to the eigenvalues $\omega_i$ of  $Z$: 
\begin{equation}
\label{eq: dual_diagonal} Z= \widetilde{Q}\operatorname{diag}(0,\dots,0,\omega_{n-s+1},\dots,\omega_{n})\widetilde{Q}^{T}.\end{equation}
\end{definition} 
\begin{definition}[Non-degeneracy]
We say that a primal-dual feasible point $(X,Z)$ is \textit{non-degenerate} if $X$ is primal non-degenerate and $Z$ is dual non-degenerate.
\newline
\end{definition}
Our interest in non-degeneracy is motivated by the following result:
\begin{proposition}
\label{thm: non_deg<->uniq}
\phantom{a}\newline
\begin{description}
    \item[\qquad1.] If $(X^*,Z^*)$ is a primal-dual non-degenerate optimal point then $(X^*,Z^*)$ is the unique \\
    \phantom{.}\qquad primal-dual optimal point for $(P,D)$.
    \item
    \item[\qquad2.] Under strict complementarity, if $(X^*,Z^*)$ is a primal-dual unique optimal point then\\ \phantom{.}\qquad$(X^*,Z^*)$ is a non-degenerate primal-dual optimal point for $(P,D)$.
    \item[]
\end{description}
\end{proposition}
\begin{remark}
For a given point $(X,Z)$, there exist linear algebraic conditions to check whether it is non-degenerate or not (see Theorems 6 and 9 in \cite{alizadeh1997complementarity}).
\end{remark}

\subsection{Set-valued analysis for {parametric} SDP}
\label{subsec: set_valued_analysis}
 
We are interested in studying the trajectories of solutions to the primal {parametric} SDP
\[
\tag{P$_t$}
\begin{aligned}
\min_{X\in\mathbb{S}^n}\ &\langle C(t), X\rangle\\
s.t. \ \ & \mathcal{A}(t)[X]=b(t)\\
&X\succeq0
\end{aligned}
\]
with time parameter $t\in T=(t_i,t_f) \subset \mathbb{R}$. 
For a given value of $t$ the dual SDP is 
\[
\label{dual_sdp}
\tag{D$_t$}
\begin{aligned}
\max_{y\in\mathbb{R}^m,\, Z\in\mathbb{S}^n}\ &\langle b(t),y\rangle\\
s.t. \quad\ \ & \mathcal{A}^*(t)[y]+Z=C(t)\\
&Z\succeq0.
\end{aligned}
\]

\begin{definition}[Set-valued maps]
A \textit{set-valued map} $F$ from a set $T$ to another set $X$ maps a point $t\in T$ to a non-empty subset of $F(t)\subseteq X$. In symbols:
\begin{align*}
    F:\,T&\rightrightarrows X\\
    t&\mapsto F(t)\subseteq X.
\end{align*}
We say that $F$ is \textit{single-valued}
at $t\in T$ if $F(t)$ is a singleton. We say that $F$ is \textit{multi-valued}
at $t\in T$ if $F(t)$ is neither empty nor a singleton. 
\\
\end{definition}

\par
Given a primal-dual pair of {parametric} SDPs $(P_t,D_t)$, we can now define the primal and dual \textit{feasible set-valued maps}:
\begin{align*}
    &\mathcal{P}(t)=\{X\in\mathbb{S}^n\;\vert\; \mathcal{A}(t)[X]=b(t),\;X\succeq0\},\\
    &\mathcal{D}(t)=\{\hspace{.025cm}Z\hspace{.025cm}\in\mathbb{S}^n\;\vert\; \mathcal{A}^*(t)[X]+Z=C(t),\;y\in\mathbb{R}^m,\;Z\succeq0\}. 
\end{align*} 
The primal and dual\textit{ optimal value functions} are defined as 
\begin{align*}
    &p^*(t)=\min_{X\in\mathbb{S}^n}\{\langle C(t), X\rangle\;\vert\;\mathcal{A}(t)[X]=b(t),\;X\succeq0\},
    \\
    &d^*(t)=\max_{Z\in\mathbb{S}^n}\{\langle b(t),y\rangle\;\vert\;\mathcal{A}^*(t)[y]+Z=C(t),\;y\in\mathbb{R}^m,\;Z\succeq0\}.
\end{align*} 
Finally, the primal and dual \textit{optimal set-valued maps} are 
\begin{align*}
    \mathcal{P}^*(t)&=\{X\in \mathcal{P}(t)\;\vert\; \langle C(t), X\rangle=p^*(t)\},\\
    \mathcal{D}^*(t)&=\{\hspace{.025cm}Z\hspace{.025cm}\in \mathcal{D}(t)\;\vert\; \langle b(t),y\rangle=d^*(t),\;\mathcal{A}^*(t)[y]+Z=C(t),\;y\in\mathbb{R}^m\}. 
\end{align*}   
Continuity properties of set-valued maps can be defined in terms of
outer and inner limits, leading to the notion of Painlevé-Kuratowski continuity. First, we introduce the notion of inner and outer limits of a set-valued map.
\newline

\par
\begin{definition}[Inner and outer limits]
Given a set-valued map $F:\,T\rightrightarrows X$, its \textit{inner limit} at $\hat t \in T$ is defined as
\[
\liminf_{t\to \hat t}F(t):= \left\{\hat x\;\big\vert\;\forall\{t_k\}_{k=1}^\infty\subseteq T\text{ such that }t_k\to \hat t,\;\exists\{x_k\}_{k=1}^\infty\subseteq X  ,\;x_k\to \hat x\text{ and } x_k\in F(t_k)\right\},
\]
while its \textit{outer limit} at $\hat t \in T$ is defined as
\[\limsup_{t\to \hat t}F(t):= \left\{\hat x\;\big\vert\;\exists\{t_k\}_{k=1}^\infty\subseteq T\text{ such that }t_k\to \hat t,\;\exists\{x_k\}_{k=1}^\infty\subseteq X  ,\;x_k\to \hat x\text{ and } x_k\in F(t_k)\right\} .\]
\end{definition}
  
\begin{definition}[Painlevé-Kuratowski continuity] 
Let $F:\,T\rightrightarrows X$ be a set-valued map. We say that $F$  is \textit{outer semi-continuous} at $\hat t\in T$ if 
\begin{equation*}
    \limsup_{t\to \hat t}F(t)\subseteq F(\hat t).
\end{equation*}
We say that $F$ is \textit{inner semi-continuous} at $\hat t\in T$ if 
\[ 
\liminf_{t\to \hat t}F(t)\supseteq F(\hat t).
\] 
Finally, we say that $F$ is \textit{Painlevé-Kuratowski continuous} at $\hat t$ if it is both outer and inner semi-continuous at $\hat t$. 

\end{definition}
\begin{remark}[Continuity]
Note that a single-valued map $F:T\to X$ is continuous in the usual
sense at a point $x\in X$ if and only if it is Painlevé-Kuratowski  continuous at  $x\in X$  as a multi-valued
map $F:T\rightrightarrows X$. Thus, without ambiguity, we will refer to Painlevé-Kuratowski continuity simply as continuity.
\\
\end{remark}

In the following, we list some continuity results on the feasible and optimal set-valued maps. The proof of Theorem \ref{thm: feasible_sets_isc} in the primal version is an original contribution of this paper.
\begin{theorem}[Example 5.8 in \cite{rockafellar2009variational}]
\label{thm: feasible_sets_osc}
If $\mathcal{A}(t)$, $b(t)$ and $C(t)$ are continuous functions of $t$ (see Assum. \ref{ass: data_continuity} in Section 3), then the feasible set-valued maps  $\mathcal{P}(t)$ and $\mathcal{D}(t)$ are outer semi-continuous at any $t\in T$. 
\end{theorem}

\begin{theorem}
\label{thm: feasible_sets_isc}
Assume that strict feasibility holds at any $t\in T$ (see Assum. \ref{ass: strict_feas} in Section 3), that the linear operator $\mathcal{A}(t)$ is surjective for every  $t\in T$, that the norm of $\mathcal{A}(t)$ and the norm of its pseudo-inverse $\mathcal{A}^*(t)\big(\mathcal{A}(t)\mathcal{A}^*(t)\big)^{-1}$ are uniformly bounded in $t$ (see Assum. \ref{ass: licq} in Section 3), and that $\mathcal{A}(t)$, $b(t)$ and $C(t)$ are continuous functions of $t$ (see Assum. \ref{ass: data_continuity} in Section 3). Then the set-valued maps $\mathcal{P}(t)$ and $\mathcal{D}(t)$ are inner semi-continuous for every $t\in T$.
\end{theorem}
\proof{Proof.}
For the dual case, we refer to Lemma 1 in \cite{hauenstein2019computing} for a version of this theorem where only the matrix $C$ depends on the parameter and this dependence is linear. We prove the primal case in the more general case where the left hand side $\mathcal{A}(t)$ is time-dependent and continuous and the right hand side $b(t)$ is continuous. The dual case can be proven in an analogous way. Fix $\hat t\in T$ and  $\hat X\in\mathcal{P}(\hat t)$. Given a sequence of times $\{t_k\}_{k=1}^\infty$ with
$t_k\to\hat t$, we will construct a convergent sequence  $X_k\to\hat X$ so that $X_k\in\mathcal{P}(t_k)$ for all
sufficiently large values of $k$. If $\hat X\succ0 $ we define 
\[
X_k:=\hat X + \mathcal{A}^*(t_k)\big(\mathcal{A}(t_k)\mathcal{A}^*(t_k)\big)^{-1}\left(b(t_k)-b(\hat t)\right).
\]
The definition is well posed because under the assumptions of the theorem the operator $\mathcal{A}(t_k)$ has full rank, thus $\mathcal{A}(t_k)\mathcal{A}^*(t_k)$ is invertible. Clearly, $\mathcal{A}(t_k)[X_k]=b(t_k)$. Furthermore, we have that $\|X_k-\hat X\|_F=\|\mathcal{A}^*(t_k)\big(\mathcal{A}(t_k)\mathcal{A}^*(t_k)\big)^{-1}\left(b(t_k)-b(\hat t)\right)\|_F\le C_{\mathcal{A}}\|b(t_k)-b(\hat t)\|\to0$ for some constant $C_{\mathcal{A}}$ (which exists by the hypothesis of uniform boundedness) and by continuity of $b(t)$, so that $X_k\to\hat X$ and $X_k\succeq0$ for sufficiently large $k$. 
If $\hat X\succeq0 $ and its smallest eigenvalue $\lambda_{\min}(\hat X)$ is zero, we define 
\[
X_k:=(1-\alpha_k)\hat X + \alpha_k\bar{X} + \mathcal{A}^*(t_k)\big(\mathcal{A}(t_k)\mathcal{A}^*(t_k)\big)^{-1}\left(b(t_k)-b(\hat t)\right) 
\]
for a fixed $\bar{X}\in\mathcal{P}(\hat t)$ such that $\bar{X}\succ0$, which exists by the strict feasibility assumption, and for a sequence $\{\alpha_k\}_{k=1}^\infty\subseteq[0,1]$ which we shall conveniently define in the following. 
 Clearly, $\mathcal{A}(t_k)[X_k]=b(t_k)$ and hence we only need to prove that $X_k\succeq0$ or, equivalently, that 
\begin{gather*} 
\lambda_{\min}\left((1-\alpha_k)\hat X + \alpha_k\bar{X} + \mathcal{A}^*(t_k)\big(\mathcal{A}(t_k)\mathcal{A}^*(t_k)\big)^{-1}\left(b(t_k)-b(\hat t)\right)\right)\ge0,
\end{gather*}
which, thanks to Weyl's inequality (see e.g Theorem 1 in \cite{franklin2012matrix}, Section 6.7) holds if
\begin{gather*}
 \alpha_k\lambda_{\min} (\bar{X}) + \lambda_{\min}\left(\mathcal{A}^*(t_k)\big(\mathcal{A}(t_k)\mathcal{A}^*(t_k)\big)^{-1}\left(b(t_k)-b(\hat t)\right)\right)
\ge0.
\end{gather*}
Rearranging:
\[
\alpha_k \ge-\frac{ \lambda_{\min}\left(\mathcal{A}^*(t_k)\big(\mathcal{A}(t_k)\mathcal{A}^*(t_k)\big)^{-1}\left(b(t_k)-b(\hat t)\right)\right)}{\lambda_{\min} (\bar{X})} .
\]
We then define $\alpha_k:=\max\{0,\beta_k\}$, where
\[
\beta_k:=-\frac{ \lambda_{\min}\left(\mathcal{A}^*(t_k)\big(\mathcal{A}(t_k)\mathcal{A}^*(t_k)\big)^{-1}\left(b(t_k)-b(\hat t)\right)\right)}{\lambda_{\min} (\bar{X})} .
\]
For sufficiently large $k$, $\beta_k\le1$, so that $\{\alpha_k\}_{k=1}^\infty\subseteq[0,1]$ and thus $X_k\in\mathcal{P}(t_k)$, since $\beta_k\to0$, $\alpha_k\to0$ and $X_k\to\hat X$.
\newline
\endproof

Theorems \ref{thm: feasible_sets_osc} and \ref{thm: feasible_sets_isc} show that the primal and dual feasible set-valued maps $\mathcal{P}(t)$ and $\mathcal{D}(t)$ are always continuous, under the assumptions of Theorem \ref{thm: feasible_sets_isc}.
Naturally, we now investigate the inner and outer semi-continuity of the optimal set-valued maps.  We have:

\begin{theorem}[Theorem 8 in \cite{hogan1973point}]
\label{thm: optimal_sets_osc}
If $\mathcal{A}(t)$, $b(t)$ and $C(t)$ are continuous functions of $t$ (see Assum. \ref{ass: data_continuity} in Section 3) and the primal and dual feasible set-valued maps are continuous, then the optimal set-valued maps  $\mathcal{P}^*(t)$ and $\mathcal{D}^*(t)$ are outer semi-continuous at any $t\in T$. 
\end{theorem}

However, in general, it is not true that the optimal set-valued maps $\mathcal{P}^*(t)$ and $\mathcal{D}^*(t)$ are inner semi-continuous. Still, the set of $t\in T$ such that $\mathcal{P}^*(t)$ or $\mathcal{D}^*(t)$ fails to be inner semi-continuous, is of \textit{first category}, i.e., countable and nowhere dense.
\begin{theorem}[Theorem 5.55 in \cite{rockafellar2009variational}]
\label{thm: loss_inner_semi_cont}
The subset of points $t\in T$ at which $\mathcal{P}^*(t)$ or $\mathcal{D}^*(t)$ fails to be inner semi-continuous (and hence continuous) is the union of countably many sets that are nowhere dense in $T$. In particular, it has empty interior. 
\end{theorem} 

However, if the optimal set is single-valued, then it is continuous everywhere. 
In order to show this, we first introduce a lemma which guarantees the local uniform boundedness of $\mathcal{P}^*$ and $\mathcal{D}^*$. 
\begin{lemma}[Lemma 3.2 in \cite{sekiguchi2021perturbation}]
\label{lem: local_uniform_boundedness}
    Assume that strict feasibility holds at any $t\in T$ (see Assum. \ref{ass: strict_feas} in Section 3) and that $\mathcal{A}(t)$, $b(t)$ and $C(t)$ are continuous functions of $t$ (see Assum. \ref{ass: data_continuity} in Section 3), then $\mathcal{P}^*(t)$ and $\mathcal{D}^*(t)$ are locally uniformly bounded at any $t\in T${, i.e., for every $t\in T$ there exists compact sets $C_p, C_d$ and $\delta\geq0$ such that $\mathcal{P}^*(s)\subseteq C_p$ and $\mathcal{D}^*(s)\subseteq C_d$ for all $s\in [t-\delta,t+\delta]$}.
\end{lemma}
\proof{Proof.}
Since we assume that primal-dual strict feasibility holds at any $t\in T$, the assumptions of both Lemmas 3.1 and 3.2 in \cite{sekiguchi2021perturbation} are satisfied at any $t\in T$.
\endproof

\begin{proposition}[Corollary 8.1 in \cite{hogan1973point}]
\label{prop: single_valued_continuous}
Assume that strict feasibility holds at any $t\in T$ (see Assum. \ref{ass: strict_feas} in Section 3), so that by Lemma \ref{lem: local_uniform_boundedness} $\mathcal{P}^*(t)$ and $\mathcal{D}^*(t)$ are locally uniformly bounded at any $t\in T$, and that $\mathcal{A}(t)$, $b(t)$ and $C(t)$ are continuous functions of $t$ (see Assum. \ref{ass: data_continuity} in Section 3). If $\mathcal{P}^*(t)$ is single-valued at $\hat{t}$, then $\mathcal{P}^*(t)$  is continuous at $\hat{t}$. The same holds for $\mathcal{D}^*(t)$.
\end{proposition} 

\subsection{Regularity properties of the {parametric} SDP optimal set-valued map}
\label{subsec: Regularity_properties_of_optimal_set_map}
Given a primal-dual pair of {parametric} SDPs $(P_t,D_t)$, we denote a primal-dual point by $(X,Z,t)$.
If at a fixed value of the parameter $\hat t\in T$ there exists a primal-dual non-degenerate optimal point $(X^*,Z^*)$, then, by Proposition \ref{thm: non_deg<->uniq}
$(X^*,Z^*)$ is a unique primal-dual optimal point, and by Proposition
\ref{prop: single_valued_continuous}, around  $\hat t$ the primal and dual optimal set-valued maps are continuous single-valued functions.
Under strict complementarity, these functions are analytic. Below we provide details of this fact.
\newline 

The optimality conditions \eqref{eq: KKT} for $(X,Z,t)$ to be a solution of $(P_t,D_t)$ at a fixed value of the parameter $t\in T$ can be equivalently written as 
\begin{equation}
\label{eq: optimality_conds}
    F(X,y,Z,t):=\begin{pmatrix} \Tilde{\mathcal{A}}(t)\operatorname{svec}(X)-b(t),
    \\
    \Tilde{\mathcal{A}}(t)^Ty+\operatorname{svec}(Z)-\operatorname{svec}(C(t)) \\
    \frac{1}{2}\operatorname{svec}\left(XZ+ZX\right)
    \end{pmatrix}=0,
\end{equation}
\begin{equation}
\label{eq: sdp_feasibility}
X,Z\succeq0
\end{equation}
for some $y\in\mathbb{R}^m$,
where $\Tilde{\mathcal{A}}(t):=\left(\operatorname{svec}(A_1(t)),\dots.\operatorname{svec}(A_m(t))\right)^T$ and $\operatorname{svec}(X)$ denotes a linear map stacking the upper triangular part of $X$, where the off-diagonal entries are multiplied by $\sqrt{2}$:
\[\operatorname{svec}(X):=\left(X_{11}, \sqrt{2} X_{12}, \ldots, \sqrt{2} X_{1n}, X_{22}, \sqrt{2} X_{23}, \ldots, \sqrt{2} X_{2 n}, \ldots, X_{n n}\right)^{T}\]
so that $\langle X,X\rangle=\operatorname{svec}(X)^T \operatorname{svec}(X)$.

\begin{definition}[Singular points]
\label{def: non_singular_points}
We say that a point $(X,y,Z)$ is \textit{singular} at $t$ if the Jacobian w.r.t $(X,y,Z)$ of $F$ at $(X,y,Z,t)$
\begin{equation}
\label{eq: jacobian} 
J_F(X,y,Z,t)=
\begin{pmatrix}
    \Tilde{\mathcal{A}}(t)&0&0
    \\
    0&\Tilde{\mathcal{A}}^T(t)&I_{\tau(n)}
    \\
    Z\otimes_sI_n&0&I_n\otimes_sX
\end{pmatrix}
\end{equation}
is not invertible, where $\otimes_s$ denotes the symmetric Kronecker product between two $n\times n$ matrices $A$ and $B$ and is defined by 
\[
(A\otimes_sB)\operatorname{svec}(H)=\frac{1}{2}(AHB^T+BHA^T)\text{\quad for any } H\in\mathbb{S}^n.
\]
Otherwise, we say that $(X,y,Z)$ is \textit{non-singular} at $t$.
 \\
 \end{definition}
\begin{definition}[Singular times]
 \label{def: non_singular_times}
 We say that a time $\hat{t}$ is \textit{singular} if there exists a singular point $(X,y,Z)$ at $\hat{t}$ such that $F(X,y,Z,\hat{t})=0$. Otherwise, we say that $\hat{t}$ is \textit{non-singular}. Furthermore, following \cite{hauenstein2019computing}, we say that a non-singular time $\hat t$ is \textit{generic} if the data tuple $(\mathcal{A}(\hat t),b(\hat t),C(\hat t))$ is generic, in the sense that the number of solutions for \eqref{eq: optimality_conds} matches the generic number of solutions (see Theorem 7 in \cite{nie2010}).
 \\
 \end{definition}
 Note that if $\hat{t}$ is non-singular, every point $(X,y,Z)$ such that $F(X,y,Z,\hat{t})=0$ is non-singular at $\hat{t}$. 
 \newline
 \par
 The following lemma gives equivalent conditions for a primal-dual optimal point $(X,Z)$ to be non-singular at $t$.
 \begin{lemma}[Theorem 3.1. in \cite{alizadeh1998primal}, Theorem 3.1. in \cite{Haeberly1998}]
 \label{lem: non_singular_jacobian}
 A primal-dual optimal point $(X,Z)$ is non-singular if and only if $(X,Z)$ is a strictly complementary and non-degenerate primal dual optimal solution.
 \end{lemma}
 Note that under strict complementarity part 2. of Proposition \ref{thm: non_deg<->uniq} holds. Therefore, the Jacobian of $F$ is non-singular at an optimal primal-dual solution $(X,Z,t)$ if and only if $(X,Z)$ is a \textit{unique} primal-dual optimal point satisfying strict complementarity. We use this result in the following theorem.
 \begin{theorem}
 \label{thm: single_valued_differentiable}
  Let $(P_t,D_t)$ be a primal-dual pair of {parametric} SDPs parametrized over a time interval $T$ such that primal-dual strict feasibility holds for any $t\in T$ (see Assum. \ref{ass: strict_feas} in Section 3) and assume that the data $\mathcal{A}(t),b(t),C(t)$ are continuously differentiable functions of $t$. Let $\Hat{t} \in T$ be a fixed value of the time parameter and  suppose that $(X^*, Z^*)$ is a unique primal-dual optimal and strictly complementary point for $(P_{\Hat{t}},D_{\hat{t}})$. Then there exists $\varepsilon > 0$ and a unique continuously differentiable mapping $(X^*(\cdot), Z^*(\cdot))$ defined on $(\hat{t}-\varepsilon,\hat{t}+\varepsilon)$ such that $(X^*(t), Z^*(t))$ is a unique and strictly complementary primal-dual optimal point to $(P_t,D_t)$ for all $t\in(\hat{t}-\varepsilon,\hat{t}+\varepsilon)$.
 \end{theorem}
 {

 \proof{Proof.}  
By Assumption \ref{ass: strict_feas} (primal-dual strict feasibility), for each $t\in T$ the pair of problems $(P_t,D_t)$ must have at least a primal-dual feasible and optimal solution which correspond to a solution of the KKT system \eqref{eq: optimality_conds}-\eqref{eq: sdp_feasibility}.
By Lemma \ref{lem: non_singular_jacobian}, the assumptions of strict complementarity and uniqueness ensure that at $\hat t$ we can apply the Implicit Function Theorem (see, e.g., Theorem 3.3.1 in \cite{krantz2002implicit}), so that there exists $\varepsilon'>0$ and a continuously differentiable curve $(X^*(\cdot),y^*(\cdot),Z^*(\cdot))$ on $(\hat{t}-\varepsilon',\hat{t}+\varepsilon')$ such that $(X^*(t),y^*(t), Z^*(t),t)$ is a solution of \eqref{eq: optimality_conds} for all $t\in(\hat{t}-\varepsilon',\hat{t}+\varepsilon')$ and $X^*(\hat{t})=X^*$, $Z^*(\hat{t})=Z^*$.
Due to the assumed strict complementarity, $\lambda_i(X^*(\hat{t}))~\hspace{-0.1cm}\cdot~\hspace{-0.1cm}\lambda_i(Z^*(\hat{t})) = 0$ and $\lambda_i(X^*(\hat{t})) +\lambda_i(Z^*(\hat{t}))>0$,  where $\lambda_i(\cdot)$ denotes the $i$-th smallest eigenvalue of a matrix, and from the continuity of the eigenvalues of $X^*(t)$ and $Z^*(t)$ with respect to $t$, the non-zero eigenvalues remain non-zero and positive for a small enough perturbation of $\hat{t}$. On~the~other hand, the last equation in \eqref{eq: optimality_conds} implies $\lambda_i(X^*(t)) \cdot\lambda_i(Z^*(t)) = 0$ for $t\in(\hat{t}-\varepsilon',\hat{t}+\varepsilon')$, so~that that the zero eigenvalues stay constant in a small open neighborhood of $\hat{t}$. In other words, the perturbed solutions remain positive semidefinite and strictly complementary, hence satisfying~\eqref{eq: sdp_feasibility}. 
Finally, by observing that the Jacobian \eqref{eq: jacobian} stay non-singular in a small open neighborhood of $\hat t$, and by Lemma \ref{lem: non_singular_jacobian}, we can conclude that $(X^*(t),Z^*(t))$ is a unique optimal solution for $(P_t,D_t)$ for~$t\in(\hat{t}-\varepsilon,\hat{t}+\varepsilon)$
and a small enough $\varepsilon'>\varepsilon>0$.
\\
}
   
By adding further assumptions, one can improve the information given by Theorem \ref{thm: loss_inner_semi_cont} on the cardinality of the singular points set and prove that the number of singular points of \eqref{eq: optimality_conds} is finite. 
\begin{theorem}[Proposition 5 in \cite{hauenstein2019computing}]
\label{thm: finiteness_of_singularity}
For the primal-dual {parametric} SDPs $(P_t,D_t)$,
assume that there exists a generic non-singular time (see Def. \ref{def: non_singular_times}) and that the data $\mathcal{A}(t),b(t),C(t)$ are polynomial functions of $t$.
Then the set of values of the time parameter $t$ at which the primal-dual optimal point is either not unique or not strictly complementary is finite.
\end{theorem}
\proof{Proof.}
Following the proof of Proposition 5 in \cite{hauenstein2019computing} we deduce that $T$ only contains a finite number of singular points, implying that the number of optimal primal-dual points that are not strictly complementary or non-unique, hence singular by Lemma \ref{lem: non_singular_jacobian}, is finite
\newline
\endproof 

Thus, under the assumption of Theorem \ref{thm: finiteness_of_singularity}, the values of $t$ at which strict complementarity or uniqueness of the primal-dual solution is lost is finite. In particular, the values of $t$ at which $\mathcal{P}^*(t)$ or $\mathcal{D}^*(t)$ fails to be inner semi-continuous (and hence fails to be continuous) are finite. It also implies that wherever $\mathcal{P}^*(t)$ defines a continuous curve of unique optima, the values of $t$ at which $\mathcal{P}^*(t)$ fails to be differentiable are finite. The same holds for $\mathcal{D}^*(t)$.

\section{A complete classification of optimal points}
The focus of our study is first put on values $t^*$ of the time parametrization interval $T$ at which strict complementarity or uniqueness of the primal-dual optimal point is lost. In other words, these are singular points preceded by non-singular points. By Theorem \ref{thm: finiteness_of_singularity} such points are finite. There, the trajectory described by the primal and dual optimal sets can exhibit a restricted number of irregular behaviors. By an irregular behavior we mean any situation that differs from the solution following a uniquely well-defined smooth curve. Describing these situations is the goal of this Section. If, instead, all primal-dual optimal points $(X,Z,t)$ are singular for every $t\in T$, the number of possible types of irregular behaviors grows. In our main Theorem \ref{thm: main_result}, we provide a complete classification of these behaviors under both cases. 
The object of our study is the trajectory of solutions to the primal SDP $(P_t)$, that is, the primal optimal set-valued map. Every result that we propose can be clearly transposed to the dual case. 
\newline

We first adopt the following standard assumptions:
\begin{assumption}[LICQ and uniform boundedness of $\mathcal{A}$]
\label{ass: licq}
The $m$ matrices $\{A_i(t)\}_{i=1,\dots,m}$ are linearly independent in $\mathbb{S}^n$ for every $t\in T$, so that the linear operator $\mathcal{A}(t)$ is surjective for every $t\in T$.
This condition is known as the \textit{linear independence constraint qualification} (LICQ). Furthermore the operator $\mathcal{A}(t)$ and its pseudo-inverse $\mathcal{A}^*(t)\big(\mathcal{A}(t)\mathcal{A}^*(t)\big)^{-1}$ have a uniformly bounded norm.
\end{assumption}
The LICQ assumption allows us to describe the dual solution just in terms of matrix $Z$ (see Remark \ref{rem: LICQ}). The assumption of uniform boundedness is needed to ensure the inner semi-continuity of the feasible set-valued maps, see Theorem \ref{thm: feasible_sets_isc}.
\begin{assumption}[Strict feasibility]
\label{ass: strict_feas}
For every $t\in T$, problem \eqref{primal_sdp} and its dual \eqref{dual_sdp} are strictly feasible.
\end{assumption} 
This assumption is standard in the SDP literature (\cite{goldfarb1999parametric}, \cite{ahmadi2021time}, \cite{hauenstein2019computing}). Strict feasibility guarantees that the primal and dual optimal sets $\mathcal{P}^*(t)$ and $\mathcal{D}^*(t)$ are non-empty and bounded for any $t\in T$(Lemma 3.2 in \cite{goldfarb1998interior}).
Checking strict feasibility of a given SDP can be done by solving another SDP and checking whether its optimal value is positive or not (see for example \cite{hauenstein2021numerical}, Theorem 3.1 and 3.5).
\begin{assumption}[Data continuity]
\label{ass: data_continuity}
Data $\mathcal{A}(t)$, $b(t)$ and $C(t)$ depend continuously on the time parameter $t$. 
\end{assumption} 

This assumption is quite general compared to those usually found in the {parametric} SDP literature, where the data are often assumed to vary linearly with respect to the time parameter. This linearity assumption is standard when one studies sensitivity properties, so that the perturbation can be assumed to be linear. Instead, our purpose is to give a geometric characterization of the points of the trajectory of solutions, in which case we can keep a high degree of generality by just assuming continuity  of the data, without any further differentiability requirement.  
\newline

Summarizing, Assumptions \ref{ass: licq}, \ref{ass: strict_feas}, and \ref{ass: data_continuity} ensure that:  
\newline
\begin{itemize}
    \item There is no duality gap: $p^*(t)= d^*(t)$ for all $t\in T$.
    \item The primal and dual optimal faces $\mathcal{P}^*(t),\, \mathcal{D}^*(t)$ are non-empty and bounded for all $t\in T$. In other words, $(P_t)$ and $(D_t)$ are both feasible and bounded.
    \item The optimal set-valued maps  are outer semi-continuous at any $t\in T$.
    \item The subset of $T$ where the optimal set-valued map fails to be inner semi-continuous has empty interior and it is the union of countably many sets that are nowhere dense in $T$.
    \item[]
\end{itemize} 
 
Equipped with the results of the previous section, we introduce a classification into six different types of primal optimal points according to the behavior of the optimal set-valued map at these points. Our purpose is to study irregularities arising after an interval where the optimal set-valued map has regular behavior. We hence classify points for which the optimal set-valued map on a left neighborhood is unique and thus continuous.
\newline
 
Let $(P_t,D_t)$ be a primal-dual pair of {parametric} SDPs with $t\in T$. For a fixed $t^*\in T$, we consider a primal optimal point $(X^*,t^*)$ for $(P_{t^*})$. Based on the behavior of the primal optimal set-valued map $\mathcal{P}^*(t)$ at $t^*$, we can distinguish between six different cases. According to these cases we classify the primal point $(X^*,t^*)$ into six different types. This can be done analogously for the dual case.
\newline

\begin{definition}[Regular point]
\label{def: type_regular}
At a \textit{regular point} $(X^*,t^*)$,  $\mathcal{P}^*(t^*)=\{X^*\}$ and there exists $\varepsilon>0$ such that
        \begin{itemize}
        \item $\mathcal{P}^*(t)$ is single-valued and continuous for every $t\in(t^*-\varepsilon,t^*+\varepsilon)$, for some $\varepsilon>0$,
        \item $\mathcal{P}^*(t)$ is differentiable at $ t^*$.\\
    \end{itemize}
\end{definition}

\begin{remark}
\label{rem: regular_point}
Note that a primal optimal point $(X^*,t^*)$ for $(P_{t^*})$ for which there exists a dual optimal point $(Z^*,t^*)$ for $(D_{t^*})$ such that $(X^*,Z^*,t^*)$ is a non-singular point for $(P_{t^*},D_{t^*})$, is necessarily a regular point. This follows directly from Theorem \ref{thm: single_valued_differentiable} and Lemma \ref{lem: non_singular_jacobian}. The converse does not hold in general.
\\
\end{remark}
\begin{definition}[Non-differentiable point]
\label{def: type_non_diff}
At a \textit{non-differentiable point} $(X^*,t^*)$, $\mathcal{P}^*(t^*)=\{X^*\}$ and there exists $\varepsilon>0$ such that 
    \begin{itemize}
        \item $\mathcal{P}^*(t)$ is single-valued and continuous for every $t\in(t^*-\varepsilon,t^*+\varepsilon)$, 
        \item $\mathcal{P}^*(t)$ is \textit{not} differentiable at $ t^*$.
        \\
    \end{itemize}   
\end{definition}

\begin{definition}[Discontinuous isolated multiple point]
\label{def: type_isol_non_uni}
At a \textit{discontinuous isolated multiple point} $(X^*,t^*)$, $X^*\in\mathcal{P}^*(t^*)$ and there exists $\varepsilon>0$ such that
    \begin{itemize}
        \item $\mathcal{P}^*(t)$ is single-valued and continuous for every $t\in(t^*-\varepsilon,t^*)\cup(t^*,t^*+\varepsilon)$, 
        \item $\mathcal{P}^*(t)$ is multi-valued at $ t^*$.
        \\
    \end{itemize}
\end{definition}

\begin{definition}[Discontinuous non-isolated multiple point]
\label{def: type_non_isol_non_uni}
At a \textit{discontinuous non-isolated multiple point} $(X^*,t^*)$, $X^*\in\mathcal{P}^*(t^*)$ and there exists $\varepsilon>0$ such that 
        \begin{itemize}
        \item $\mathcal{P}^*(t)$ is continuous at any $t\in(t^*-\varepsilon,t^*)\cup(t^*,t^*+\varepsilon)$,
        \item $\mathcal{P}^*(t)$ is single-valued for every $t\in(t^*-\varepsilon,t^*)$,
        \item $\mathcal{P}^*(t)$ is multi-valued for every $t\in[t^*,t^*+\varepsilon)$.
        \\
    \end{itemize}  
\end{definition}
 
\begin{remark}
Let $(X^*_1,t_1^*)$ be a discontinuous isolated multiple point and $(X^*_2,t_2^*)$ a discontinuous non-isolated multiple point. Then by definition the optimal solution is not unique neither at $t^*_1$ nor at $t^*_2$. Thus, a loss of inner semi-continuity of the optimal set-valued map $\mathcal{P}^*(t)$ must occur both at $t^*_1$ and at $t^*_2$. However, while for any $\varepsilon>0$ the set of points $t\in (t^*_2-\varepsilon,t^*_2+\varepsilon)$ where the optimal set $\mathcal{P}^*(t)$ is multi-valued has a non-empty interior, there always exists a $\bar\varepsilon>0$ such that the set of points $t\in (t^*_1-\bar\varepsilon,t^*_1+\bar\varepsilon)$ where the optimal set $\mathcal{P}^*(t)$ is multi-valued has empty interior. This observation suggests the choice of the terms ``isolated" and ``non-isolated".
\\
\end{remark}

\begin{definition}[Continuous bifurcation point]
\label{def: type_bif}
At a \textit{continuous bifurcation point} $(X^*,t^*)$, $\mathcal{P}^*(t^*)=\{X^*\}$ and
there exists $\varepsilon>0$ such that

        \begin{itemize}
        \item $\mathcal{P}^*(t)$ is continuous at any $t\in(t^*-\varepsilon,t^*+\varepsilon)$,
        \item $\mathcal{P}^*(t)$ is single-valued for every $t\in(t^*-\varepsilon,t^*]$,
        \item $\mathcal{P}^*(t)$ is multi-valued for every $t\in(t^*,t^*+\varepsilon)$.
    \end{itemize}
     In particular, there exist at least two distinct continuous curves 
     \begin{align*}  
    \begin{array}{rcl}
         X_1:(t^*,t^*+\varepsilon)&\to&\mathbb{S}^n  \\
         t& \mapsto& X_1(t)
    \end{array} 
    \begin{array}{rcl}
     X_2:(t^*,t^*+\varepsilon)&\to&\mathbb{S}^n  \\
         t& \mapsto& X_2(t)
    \end{array} 
    \end{align*}
     such that  $X_1(t)$ and $X_2(t)$ are two distinct points of $\mathcal{P}^*(t)$ for every $t \in(t^*,t^*+\varepsilon)$ 
     and $\lim_{t\to t^{*+}}X_1(t)=\substack{}\lim_{t\to t^{*+}}X_2(t)=X^*$. In this sense, a continuous bifurcation point can be thought as a continuous loss of uniqueness from a single branch into two or more branches.
     \\
\end{definition}

\begin{definition}[Irregular accumulation point]
\label{def: type_accu}
At an \textit{irregular accumulation point} $(X^*,t^*)$, $X^*\in\mathcal{P}^*(t^*)$ and there exists $\varepsilon>0$ such that
    \begin{center}
          \begin{itemize}
        \item $\mathcal{P}^*(t)$ is single-valued and continuous for every $t\in(t^*-\varepsilon,t^*)$ 
    \end{itemize}
    \end{center}
    and for any $\delta>0$ at least one of the following is true:
    \begin{itemize}  
        \item there exists a sequence of times $\{t_k\}_{k=1}^\infty\subseteq(t^*,t^*+\delta)$ at which a loss of inner semi-continuity occurs and  $\lim_{k\to\infty}t_k=t^*$. At these times, either a discontinuous isolated multiple point or a discontinuous non-isolated multiple point appears.
        \item there exists a sequence of times $\{t_k\}_{k=1}^\infty\subseteq(t^*,t^*+\delta)$ at which a continuous bifurcation occurs and $\lim_{k\to\infty}t_k=t^*$.
        \\
    \end{itemize}
\end{definition}
When convenient, instead of saying that $(X^*,t^*)$ is a regular point, we will say that $X^*$ is a regular point at $t^*$. The same applies to all the other types of points that we defined.
 \\
\begin{remark}
\label{rem: reverse_time}
The above definitions consider points whose sufficiently small \textit{left} time neighborhood consists of all regular points. By a change of sign of the parameter, the definition clearly extends to points whose sufficiently small \textit{right} time neighborhood consists of all regular points.
\\
\end{remark}
\begin{remark}[Existence of a continuous selection]
The optimal set-valued map is continuous in a neighborhood of a regular, non-differentiable, or a continuous bifurcation point.
Instead, at a discontinuous isolated or non-isolated multiple point (Definitions \ref{def: type_isol_non_uni} and \ref{def: type_non_isol_non_uni}), a loss of inner semi-continuity occurs. For such points $(X^*,t^*)$ it holds $\liminf_{t\to t^{*-}}\mathcal{P}^*(t)\neq\mathcal{P}^*(t^*)$. However, in both cases, clearly only one of the following is true:
 \begin{align*}
        &(A) \lim_{t\to t^{*+}}\mathcal{P}^*(t)=\mathcal{P}^*(t^*), \\
        &(B) \liminf_{t\to t^{*+}}\mathcal{P}^*(t)\neq\mathcal{P}^*(t^*). 
    \end{align*}    
In case (A), one can select a continuous curve  $ (t^*-\varepsilon,t^*+\varepsilon)\ni t\mapsto X(t)\in\mathbb{S}^n $ such that $X(t)\in\mathcal{P}^*(t)$ for every $t \in(t^*-\varepsilon,t^*+\varepsilon)$, while in case (B) such a curve does not exist. Furthermore, for a discontinuous isolated multiple point under case (A), such a curve is unique. Also note that in case (A) it might be impossible to select a curve that is differentiable at $t^*$.
\\
\end{remark}

\begin{remark}[Comparison with \cite{guddat1990parametric}]
\label{rem: non_degeneracy}
The definition of the six different types of points was inspired by \cite[Chapter 2]{guddat1990parametric}, where a classification of solutions to univariate parametric non-linear constrained optimization problems was proposed. There, critical primal-dual points satisfying first-order optimality (or KKT) conditions for a given parametric non-linear optimization problem are classified. These points are defined as \textit{non-degenerate} if strict complementarity holds as well as the invertibility of the Hessian of the Lagrangian of the considered problem restricted to the tangent space at the point. 
We remark that this notion of non-degeneracy does not coincide with that of primal and dual non-degeneracy defined in Definitions \ref{def: primal_non_deg} and \ref{def: dual_non_deg}.
However, one can still identify an algebraic resemblance between primal non-degeneracy as defined in \ref{def: primal_non_deg} and the non-singularity of the Hessian of the Lagrangian.

In the terminology that we used, the notion of \textit{non-degeneracy} adopted by Jongen in \cite{guddat1990parametric} is analogous to \textit{non-singularity}, as defined in Definition \ref{def: non_singular_points}, as they both guarantee the applicability of the implicit function theorem, hence ensuring a regular behavior (Theorem 2.4.2 in \cite{guddat1990parametric}). Around these points the optimal set can be  parametrized by means of a single parameter and the parameterization is a differentiable map.
If a critical point is instead \textit{degenerate} then, according to which algebraic condition is not satisfied by such points, these are classified in four different types. Instead, we classified irregular points according to the behavior of the trajectory of solutions at the point considered, focusing at the possible local topological structure of points   
\end{remark}
 
\begin{theorem}[Main result]
\label{thm: main_result}
For a primal-dual pair of {parametric} SDPs $(P_t,D_t)$, let Assumptions \ref{ass: licq}, \ref{ass: strict_feas}, and \ref{ass: data_continuity} hold and consider a time $t^*\in T$, $X^*\in\mathcal{P}^*(t^*)$. If $\mathcal{P}^*(t)$ is unique for every $t\in(t^*-\varepsilon',t^*)$ for some $\varepsilon'>0$, then $(X^*,t^*)$  must be a point of a type defined in Definitions \ref{def: type_regular}, \ref{def: type_non_diff}, \ref{def: type_isol_non_uni}, \ref{def: type_non_isol_non_uni}, \ref{def: type_bif}, or \ref{def: type_accu}. The same holds for $\mathcal{D}^*(t)$
\end{theorem}

 \proof{Proof.}
 First, let $t^*\in T$ and $X^*\in\mathcal{P}^*(t^*)$.
 By hypothesis, there exists $\varepsilon'>0$ such that $\mathcal{P}^*(t)$ is single-valued and hence, by Proposition \ref{prop: single_valued_continuous}, continuous for every $t\in(t^*-\varepsilon',t^*)$. Let us perform a first binary case partition:
 \\
 \begin{description}
     \item[A] $\mathcal{P}^*(t^*)$ is a single-valued (and thus equal to $\{X^*\}$).
     \item[B] $\mathcal{P}^*(t^*)$ is multi-valued.
     \\
 \end{description} 
  Then, we  also define a three-way case partition, independent from the previous one: \\
 
 \begin{description}
     \item[1] there exists $\varepsilon''>0$ such that $\mathcal{P}^*(t)$ is single-valued for every $t\in(t^*,t^*+\varepsilon'')$.
     \item[2] there exists $\varepsilon''>0$ $\mathcal{P}^*(t)$ is multi-valued for every $t\in(t^*,t^*+\varepsilon'')$.
     \item[3] for every $\delta>0$ there exists $t',t''\in(t^*,t^*+\delta)$ such that $\mathcal{P}^*(t')$ is single-valued and $\mathcal{P}^*(t'')$ is multi-valued.\\
 \end{description}
  Combining the two partitions, we obtain one consisting of six cases:
 \\
 \begin{description}
     \item[A1] in this case $\mathcal{P}^*(t^)$ is a single-valued function defined in $(t^*-\varepsilon,t^*+\varepsilon)$, where $\varepsilon:=\min\{\varepsilon',\varepsilon''\}$, which is hence continuous by Proposition \ref{prop: single_valued_continuous}. According to whether $\mathcal{P}^*(t)$ is differentiable at $t^*$ or not, $(X^*,t^*)$ is a regular point or a non-differentiable point.
     \item[] 
     \item[A2] if there exists $\varepsilon''>0$ such that $\mathcal{P}^*(t)$ is continuous at any $t\in(t^*-\varepsilon',t^*+\varepsilon'')$ then by definition $(X^*,t^*)$ is a continuous bifurcation point (Definition \ref{def: type_bif}). Otherwise, for every $k\in\mathbb{N}$ there must exist a point $t_k\in(t^*,t^*+\frac{1}{k})$ such that a loss of inner semi-continuity occurs a $t_k$. Hence, $(X^*,t^*)$ is an irregular accumulation point (Definition \ref{def: type_accu}).
     \item[]
     \item[A3] if there exists $\varepsilon''>0$ such that $\mathcal{P}^*(t)$ is continuous at any $t\in(t^*-\varepsilon',t^*+\varepsilon'')$ then, as for any $\delta>0$ a continuous switch from unique to non-unique solutions must occur, we can construct a sequence of times $\{t_k\}_{k=1}^\infty$ at which a continuous bifurcation occurs converging to $t^*$. Otherwise, we can proceed as in case \textbf{A2} and construct a sequence of times at which a loss of inner semi-continuity occurs converging to $t^*$. Hence, $(X^*,t^*)$ is an 
     irregular accumulation point.
     \item[]
     \item[B1] in this case, simply by definition, $(X^*,t^*)$ is a discontinuous isolated multiple point (Definition \ref{def: type_isol_non_uni}).
     \item[]
     \item[B2] if there exists $\varepsilon''>0$ such that $\mathcal{P}^*(t)$ is continuous at any $t \in(t^*+\varepsilon'')$, by definition $(X^*,t^*)$ is a discontinuous non-isolated multiple point (type \ref{def: type_isol_non_uni}). Otherwise, for every $k\in\mathbb{N}$ there exists a point $t_k\in(t^*,t^*+\frac{1}{k})$ such that a loss of inner semi-continuity occurs a $t_k$. Hence, $(X^*,t^*)$ is an irregular accumulation point.
     \item[]
     \item[B3] the same discussion as in \textbf{A3}, $(X^*,t^*)$ is hence an irregular accumulation point. $\square$
     \\
 \end{description}
\endproof 
\begin{theorem} 
\label{thm: finite_bad_result}
For a primal-dual pair of {parametric} SDPs $(P_t,D_t)$, let Assumptions \ref{ass: licq}, \ref{ass: strict_feas}, and \ref{ass: data_continuity} hold. Suppose that  there exists a generic non-singular time (cf. Def. \ref{def: non_singular_times}) and that the data of $(P_t,D_t)$ are polynomial functions of $t$. Then, along the parametrization interval $T$ the number of points in times at which there is a non-differentiable point (cf. Def. \ref{def: type_non_diff}) or a discontinuous isolated multiple point (cf. Def. \ref{def: type_isol_non_uni}) for $\mathcal{P}^*(t)$ or $\mathcal{D}^*(t)$ is finite. All the other points are regular points (cf. Def. \ref{def: type_regular}) for both $\mathcal{P}^*(t)$ and $\mathcal{D}^*(t)$. Furthermore, the number of regular points where $\mathcal{P}^*(t)$ or $\mathcal{D}^*(t)$ is not continuously differentiable is finite.
\end{theorem}
 \proof{Proof.}
By Theorem \ref{thm: finiteness_of_singularity}, the hypothesis implies that the number of values of $t\in T$ at which there exists an optimal primal-dual singular point for \eqref{eq: optimality_conds} is finite. Let $S$ denote the set of such values. First, let $t^*_{ns}\in T\setminus S$. Then there exists an optimal primal-dual non-singular point $(X^*_{ns},Z^*_{ns},t^*_{ns})$. By Theorem \ref{thm: single_valued_differentiable}, both $(X^*_{ns},t^*_{ns})$ and $(Z^*_{ns},t^*_{ns})$ are regular points (cf. Def. \ref{def: type_regular} and Rem. \ref{rem: regular_point}) where both $\mathcal{P}^*(t)$ and $\mathcal{D}^*(t)$ are continuously differentiable. Now consider $t^*_s\in S$. Then there exists an optimal primal-dual singular point $(X^*_s,Z^*_s,t^*_s)$.  If at $t^*_s$ a loss of inner semi-continuity for $\mathcal{P}^*$ occurs then $\mathcal{P}^*(t^*_s)$ is multi-valued, hence $(X^*_s,t^*_s)$ is a discontinuous isolated multiple point (cf. Def. \ref{def: type_isol_non_uni}). The same holds in the dual version for $\mathcal{D}^*$ and $(Z^*_s,t^*_s)$. If instead at $t^*_s$ continuity of $\mathcal{P}^*$  is preserved, then $\mathcal{P}^*(t^*_s)$ is a singleton. According to whether $\mathcal{P}^*$ is differentiable at $t^*_s$ or not, $(X^*_s,t^*_s)$ is a regular point or a non-differentiable point (cf. Def. \ref{def: type_non_diff}). At regular points in $S$ that are differentiable, the derivative of $\mathcal{P}^*(t)$ and $\mathcal{D}^*(t)$ might yet fail to be continuous. Being in $S$, such points are in a finite number, hence proving the last sentence of the theorem. Since $\mathcal{P}^*(t^*_s)$ is a singleton, a loss of differentiability only happens when $t^*_s$ is in $S$; that is, when either $\mathcal{D}^*(t^*_s)$  is multi-valued or strict complementarity between $X^*_s$ and $Z^*_s$ fails (this follows from Lemma \ref{lem: non_singular_jacobian}). The same holds in the dual version for $\mathcal{D}^*$ and $(Z^*_s,t^*_s)$. 
\newline
\endproof 
To prove that any type of point that we  defined can actually appear, in the following section we exhibit an example of each type.
 \section{Examples}
 \noindent
 \subsection{Regular, non-differentiable and discontinuous isolated multiple points}
For $t\in T=(-3,2)$, consider the primal SDP
\begin{equation*}
    \tag{P$^1_t$}
    \begin{aligned}
        \operatorname{min} \ &\ \ tx+ty+z\\
         s.t.\  &\begin{pmatrix}
    1&&&x&&&y\\
    x&&&1&&&z\\
    y&&&z&&&1
    \end{pmatrix}\succeq0.
    \end{aligned}
\end{equation*}
 
The feasible region is known as Cayley spectrahedron. We have: \begin{equation*}
\mathcal{P}^*(t)=
\begin{cases}
 \begin{pmatrix}
    1&&&1&&&1\\
    1&&&1&&&1\\
    1&&&1&&&1
    \end{pmatrix} &\text{for $t\in(-3,-2]$},
    \\
    \\
     \begin{pmatrix}
    1&-t/2&-t/2\\
    -t/2&1&\frac{t^2}{2}-1\\
    -t/2&\frac{t^2}{2}-1&1
    \end{pmatrix} &\text{for $t\in(-2,2)\setminus\{0\}$},
    \\
    \\
    \left\{\begin{pmatrix}
    1&&a&&b\\
    a&&1&&-1\\
    b&&-1&&1
    \end{pmatrix}\Bigg|\ \begin{array}{cc}
         a+b=0  \\
         a,b\in[-1,1]
    \end{array} \right\}&\text{at $t=0$}.
\end{cases}   \end{equation*}
\begin{figure}[ht]
        \centering
        \includegraphics[width=6.5cm]{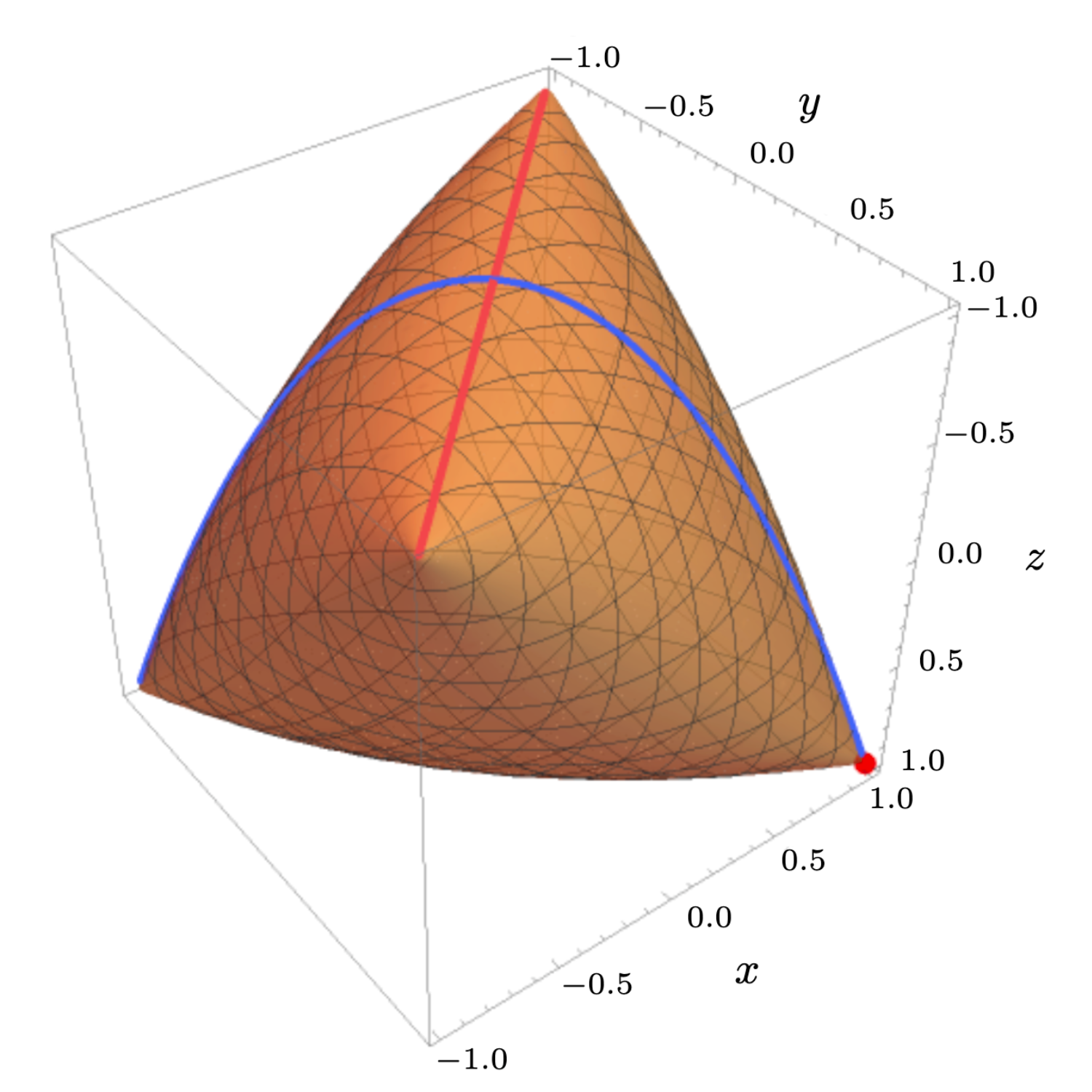}
        \caption{
        Trajectory of solutions of (P$^1_t)$. Its feasible set is time-invariant and it is the Cayley spectrahedron (orange). Its optimal set-valued map coincides with the red dot at $(1,1,1)$ for $t\in(-3,-2]$, moves along the blue curve $(-t/2,-t/2,t^2/2-1
        )$ for $t\in(-2,2)\setminus\{0\}$, and covers the whole red top 
        edge $\{(x,y,-1)|x+y=0\}$ at $t=0$.}
        \label{fig: example_1}
\end{figure}
  \noindent 
  \newline
In $(-3,-2)$, the trajectory is constant. All points are hence regular (Def. \ref{def: type_regular}). In both intervals $(-2,0)$ and $(0,2)$, the solution to (P$^1_t)$ is unique and the trajectory describes a parabolic differentiable curve  and hence all its points are also regular.
  \newline
   \noindent
   Instead, $t=-2$ is a non-differentiable point (Def. \ref{def: type_non_diff}). Indeed:
  \[
  \frac{d}{dt}\mathcal{P}^*(t)|_{t=-2^{-}}=\begin{pmatrix}
    0&0&0\\
    0&0&0\\
    0&0&0
    \end{pmatrix}\neq\begin{pmatrix}
    0&-0.5&-0.5\\
    -0.5&0&-2\\
    -0.5&-2&0 
    \end{pmatrix}=\frac{d}{dt}\mathcal{P}^*(t)|_{t=-2^{+}}.
  \]
  
  \noindent
  Moreover, at $t=0$ there is a loss of uniqueness, as $\mathcal{P}^*(0)$ is a one-dimensional face of Cayley spectrahedron. Thus, $t=0$ is a discontinuous isolated multiple point (Def. \ref{def: type_isol_non_uni}), as uniqueness is holding before for $t\in(-2,0)$ and after for $t\in(0,3)$.   
  \newline
  Consider now the SDP dual to (P$^1_t)$
   \begin{gather}
    \operatorname{max} \ \alpha+\beta+\gamma\nonumber\\
    s.t.\quad\begin{pmatrix}
    -\alpha&t/2&t/2\\
    t/2&-\beta&1/2\\
    t/2&1/2&-\gamma
    \end{pmatrix}\succeq0\tag{D$^1_t$}.
\end{gather}
  The optimal set-valued map for (D$^1_t)$ is
  \begin{equation*}
  \mathcal{D}^*(t)=
      \begin{cases}
       \begin{pmatrix}
    -t&t/2&t/2\\
    t/2&-(t+1)/2&1/2\\
    t/2&1/2&-(t+1)/2 
    \end{pmatrix}&\text{for $t\in(-3,-2)$},
    \\
    \\
     \begin{pmatrix}
    t^2/2&t/2&t/2\\
    t/2&1/2&1/2\\
    t/2&1/2&1/2
    \end{pmatrix}&\text{for $t\in[-2,2)$}.
      \end{cases}
  \end{equation*} 
  \noindent At $t=-2$, $\mathcal{D}^*(t)$ has a non-differentiable point (Def. \ref{def: type_non_diff}) too. Indeed:
  \[
  \frac{d}{dt}\mathcal{D}^*(t)|_{t=-2^{-}}=\begin{pmatrix}
    -1&0.5&0.5\\
    0.5&-0.5&0\\
    0.5&0&-0.5 
    \end{pmatrix}\neq\begin{pmatrix}
    -2&0.5&0.5\\
    0.5&0&0\\
    0.5&0&0
    \end{pmatrix}=\frac{d}{dt}\mathcal{D}^*(t)|_{t=-2^{+}}.
  \]
  
   \noindent 
  For $t\in(-3,2)\setminus\{-2\}$ the primal-dual pair of solutions is strictly complementary. Being both unique solutions for every $t\in(-3,2)\setminus\{0\}$, we conclude by Lemma \ref{lem: non_singular_jacobian} and Theorem \ref{thm: single_valued_differentiable} that for $t\in(-3,2)\setminus\{-2,0\}$ the primal-dual trajectory of solutions consists of regular points.
  \newline
  
  \noindent
  Notice that $-2$ and $0$ are singular times for the parameterization interval $T=(-3,2)$. Indeed, at $t=-2$ there is a loss of strict complementarity (the rank of both primal and dual solution is 1), while at $t=0$ there is a loss of primal uniqueness, hence a dual degenerate solution.
  \newline
  
  \noindent
  Note that this example illustrates Theorem \ref{thm: finite_bad_result}, as there exists a non-singular time $\hat{t}\in(-3,2)$ (Def. \ref{def: non_singular_times}). Take for example $\hat{t}=1$: equation \eqref{eq: optimality_conds} has a finite set of 8 solutions, which can be described as the intersections in $\mathbb{R}^6$ of 3 sets, each of which is the union of 2 hyperplanes, with 3 hyperplanes. If we set
  \begin{equation*}
   (X,Z)=\left(\begin{pmatrix}
    1&&&x&&&y\\
    x&&&1&&&z\\
    y&&&z&&&1
    \end{pmatrix},
    \begin{pmatrix}
    -\alpha&1/2&1/2\\
    1/2&-\beta&1/2\\
    1/2&1/2&-\gamma
    \end{pmatrix}\right),
  \end{equation*}
then equation \eqref{eq: optimality_conds} can be rewritten as:
\begin{equation}
    \begin{cases}
        x=\alpha+\beta-\gamma\\
        y=\alpha-\beta+\gamma\\
        z=-\alpha+\beta+\gamma\\ 
        (1+\alpha-\beta-\gamma)(1+\beta+\gamma)=0\\
        (1-\alpha+\beta-\gamma)(1+\alpha+\gamma)=0\\
        (1-\alpha-\beta+\gamma)(1+\alpha+\beta)=0.\\
    \end{cases}   
\end{equation}
 The solutions of this system are:
  \begin{equation*}
  \begin{aligned}
      &(-\tfrac{1}{2},-\tfrac{1}{2},-\tfrac{1}{2},-\tfrac{1}{2},-\tfrac{1}{2},-\tfrac{1}{2}),
      &&(1,1,1,1,1,1),
      \\
      &(1,1,-2,1,-\tfrac{1}{2},-\tfrac{1}{2}),
      && (-1,-1,1,-1,0,0),
      \\
      &(1,-2,1,-\tfrac{1}{2},1,-\tfrac{1}{2}),
      &&(-1,1,-1,0,-1,0),
      \\
      &(-2,1,1,-\tfrac{1}{2},-\tfrac{1}{2},1),
      &&(1,-1,-1,0,0,-1).
  \end{aligned}
  \end{equation*}
It is then possible to check that each of these 8 points makes the Jacobian \eqref{eq: jacobian} invertible, hence guaranteeing that $\hat{t}=1$ is a non-singular time, so that the hypothesis of Theorem \ref{thm: finite_bad_result} are satisfied. Notice that the first solution above corresponds to the optimal primal-dual solution at $\hat{t}=1$.
\newline
 
\subsection{Discontinuous non-isolated multiple points}
 For $t\in T=(-2,1)$, consider the SDP
 \begin{gather}
    \operatorname{min} \ tx+ty+z\nonumber\\
    s.t.\quad\begin{pmatrix}
    1&x&y&0\\
    x&1&z&0\\
    y&z&1&0\\
    0&0&0&1+x+y+z
    \end{pmatrix}\succeq0\tag{P$^2_t$}
\end{gather}
for which
\begin{equation*}
\mathcal{P}^*(t)=
    \begin{cases}
    \begin{pmatrix}
    1&-t/2&-t/2&0\\
    -t/2&1&\frac{t^2}{2}-1&0\\
    -t/2&\frac{t^2}{2}-1&1&0\\
    0&0&0&\frac{t^2}{2}-t
    \end{pmatrix} &\text{for $t\in(-2,0)$},
    \\
    \\
    \left\{\begin{pmatrix}
    1&a&b&0\\
    a&1&-1&0\\
    b&-1&1&0\\
    0&0&0&0
    \end{pmatrix}\Big|\ \begin{array}{cc}
         a+b=0  \\
         a,b\in[-1,1]
    \end{array} \right\}&\text{for $t\in[0,1)$}.
    
    \end{cases}
\end{equation*} 
 \begin{figure}[ht]
        \centering
        \includegraphics[width=6.5cm]{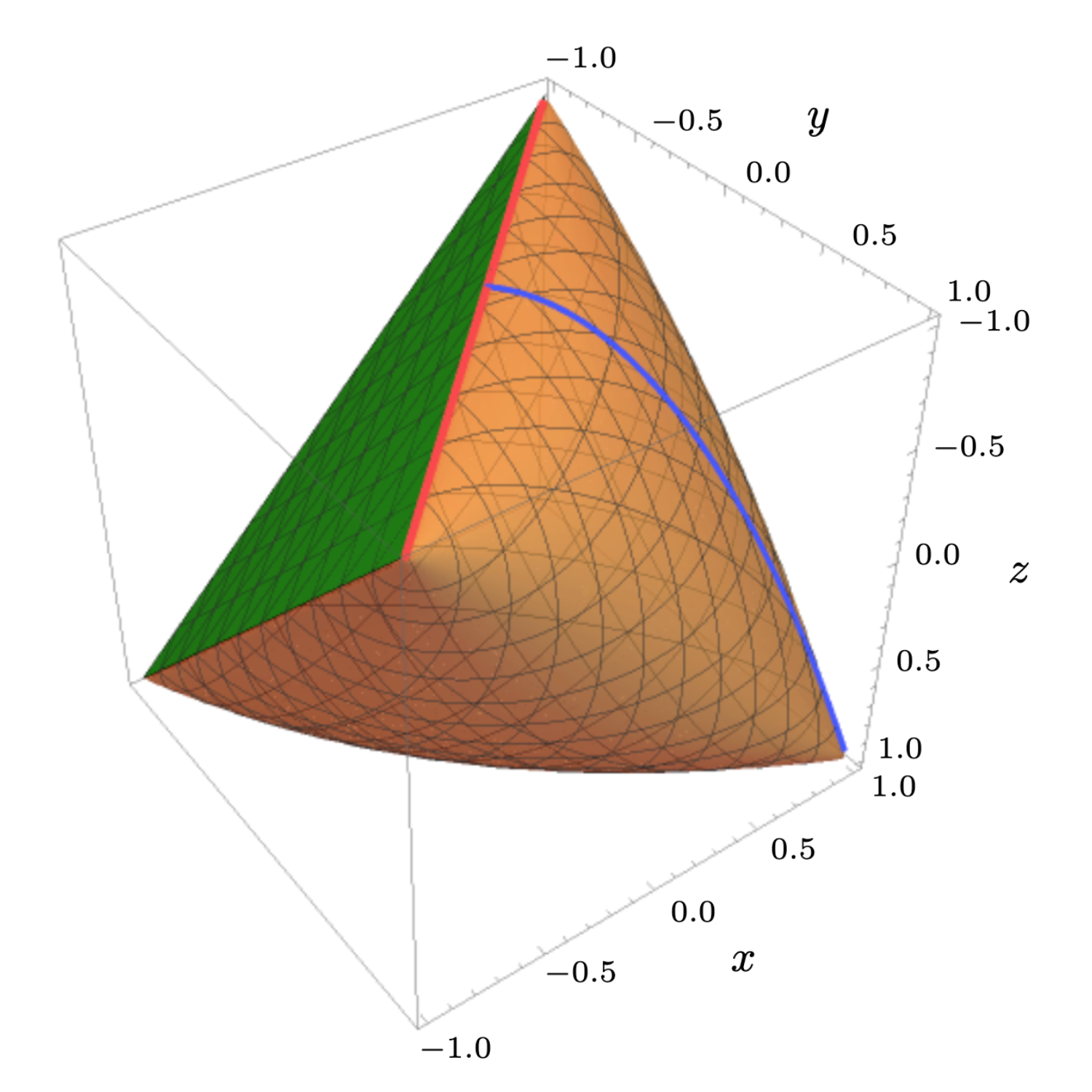}
        \caption{Trajectory of solutions of (P$^2_t)$. Its feasible set is time-invariant and it is the Cayley spectrahedron (orange) intersected with half space $\{(x,y,z)|1+x+y+z\ge0\}$ (green). Its optimal set-valued map moves along the blue curve $(-t/2,-t/2,t^2/2-1
        )$ for $t\in(-1,0)$, and covers the whole red top 
        edge $\{(x,y,z)|x+y=0,z=-1\}$ for $t\in[0,1)$.}
        \label{fig: example_2}
\end{figure}
\noindent
The optimal set-valued map $\mathcal{P}^*(t)$ is continuous for every $t\in(-2,1)\setminus\{0\}$, it is single-valued for every $t\in(-2,0)$, and it is multi-valued for every $t\in[0,1)$, as for every $t\in[0,1)$ the optimal face at $t$ is 1-dimensional. A loss of inner semincontinuity occurs at $t=0$. Hence, $t=0$ is a discontinuous non-isolated multiple point, according to Def. \ref{def: type_non_isol_non_uni}.
    
 \subsection{Continuous bifurcation point}
 For $t\in T=(-1,1)$, consider the primal SDP
 \begin{equation}
 \tag{P$^3_t$}
      \begin{aligned}
    \operatorname{min} \ &x_{11}\nonumber\\
    s.t.\ &x_{44}-x_{33}=0\\
    &x_{22}=1\\
    &2x_{12}+x_{33}+x_{44}=-t\\
    &X\succeq0
\end{aligned}
 \end{equation}
for which
\begin{equation*}
    \mathcal{P}^*(t)=
    \begin{cases}
    \left\{ \begin{pmatrix}
    0&0&0&0\\
    0&1&a&b\\
    0&a&-t/2&c\\
    0&b&c&-t/2
    \end{pmatrix}\Bigg|\ \begin{array}{cc}
         a^2+b^2+c^2\le\frac{t^2}{4}-t\\
         -\frac{t}{2}(a^2+b^2)+c^2-2abc\le\frac{t^2}{4}
    \end{array} \right\}&\text{for $t\in(-1,0)$},
    \\
    \\
    \begin{pmatrix}
    t^2/4&-t/2&0&0\\
    -t/2&1&0&0\\
    0&0&0&0\\
    0&0&0&0
    \end{pmatrix}&\text{for $t\in[0,1)$}.
    \end{cases}
\end{equation*} 
The optimal set-valued map $\mathcal{P}^*(t)$ is continuous for every $t\in(-1,1)$, it is multi-valued for every $t\in(-1,0)$, being there a 3-dimensional face, and it is single-valued for every $t\in[0,1)$. Hence $t=0$ is a continuous bifurcation point for (P$^3_t)$ according to Def. \ref{def: type_bif} (with reversed time, see Remark \ref{rem: reverse_time}).
\newline

\noindent
When there exists a continuous bifurcation point it is necessary that all the times of the parameterization interval are singular according to Def. \ref{def: non_singular_times}. In other words, at any time $t\in(-1,1)$ there exists a primal-dual point which is either degenerate or not strictly complementary. 
Indeed, the dual SDP to (P$^3_t)$ is 
\begin{gather}
\operatorname{max} \ y-tz\nonumber\\
s.t.\quad\begin{pmatrix}
1&-z&0&0\\
-z&-y&0&0\\
0&0&-x-z&0\\
0&0&0&x-z
\end{pmatrix}\succeq0\tag{D$^3_t$},
\end{gather}
which is equivalent to $\max\{y+tz\;|\;y+z^2\le0,\; -z\le x\le z\}$
and for which
 \begin{equation*}
     \mathcal{D}^*(t)=
     \begin{cases}

    \begin{pmatrix}
    1&0&0&0\\
    0&0&0&0\\
    0&0&0&0\\
    0&0&0&0
    \end{pmatrix}&\text{for $t\in(-1,0]$},
    \\
    \\
    \left\{\begin{pmatrix}
    1&t/2&0&0\\
    t/2&t^2/4&0&0\\
    0&0&-a&0\\
    0&0&0&a+t
    \end{pmatrix}\Bigg|\ a\in[-t,0] \right\}&\text{for $t\in(0,1)$}.
     \end{cases}
 \end{equation*} 
  
The dual optimal set-valued map $\mathcal{D}^*(t)$ is continuous for every $t\in(-1,1)$, single-valued for every $t\in(-1,0]$, and it is multi-valued for every $t\in(0,1)$, being there a 1-dimensional face. Thus, $t=0$ is a continuous bifurcation point for (D$^3_t)$, according to Def. \ref{def: type_bif}.
\newline

\noindent In particular, a pair of primal-dual solutions for (P$^3_t$,D$^3_t)$ is not unique, hence degenerate, for every $t\in(-1,1)\setminus\{0\}$; for $t=0$, there is a unique pair of primal-dual solutions for which however strict complementarity does not hold. This implies that all $t\in (-1,1)$ are singular times.

\begin{figure}[ht]
\centering \includegraphics[width=6.5cm]{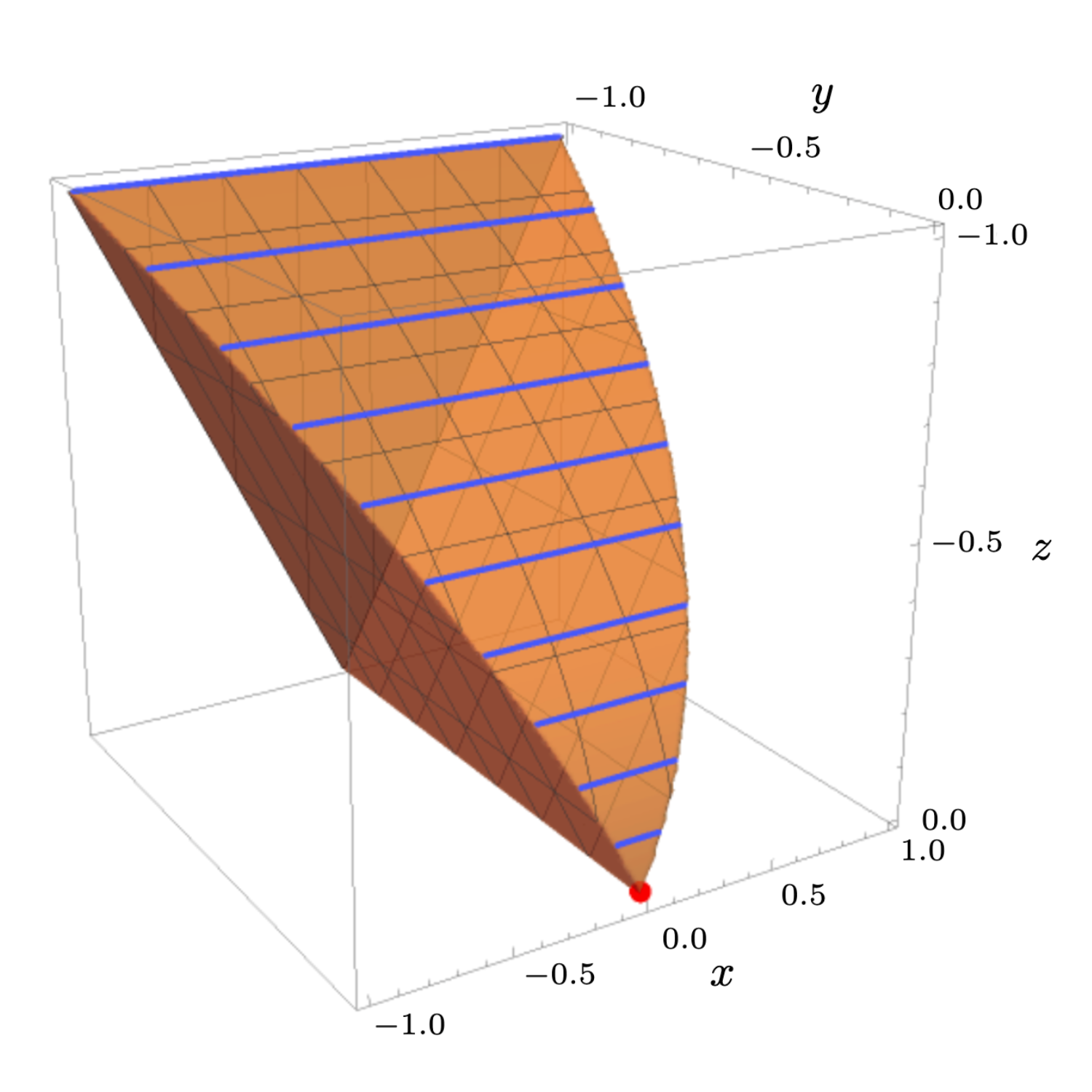} 
\caption{Trajectory of solutions of (D$^3_t)$. Its feasible set is time-invariant and it is the set  $\{(x,y,z)|y+z^2\le0,\; -z\le x\le z\}$ (orange). Its optimal set-valued map coincides with the red dot at $(0,0,0)$ for $t\in(-1,0]$. At $t=0$, $(0,0,0)$ is a continuous bifurcation point, as for every $t\in(0,1)$ the solution is multi-valued and equal to the set $\{(x,y,z)|x\in[-t/2,t/2], y=-t^2/4, z=-t/2\}$. In the picture, the blue segments illustrate the optimal multiple-valued solution for  $t=\{0.1,0.2,\dots,0.9,1\}$}
\label{fig: example_3}
\end{figure}
\subsection{Irregular accumulation points}
For $t\in T=(-1,1)$, consider the SDP
 \begin{gather}
    \operatorname{min} \ f(t)(x-y)+z\nonumber\\
    s.t.\quad\begin{pmatrix}
    1&x&y&0&0\\
    x&1&z&0&0\\
    y&z&1&0&0\\
    0&0&0&g(t)&x-y\\
    0&0&0&x-y&g(t)
    \end{pmatrix}\succeq0\tag{P$^4_t$}
\end{gather}
where 
\[
f(t):=
\begin{cases}
    t\sin{\frac{\pi}{t}}&\text{if }t>0,
    \\
    0&\text{otherwise},
\end{cases}
    \text{ and }g(t):=
\begin{cases}
    2t&\text{if }t>0,
    \\
    0&\text{otherwise}.
\end{cases}
\]
For $t\le0$ the feasible region is the intersection between Cayley spectrahedron and the plane $x-y=0$. 
For $t>0$ the feasible region is the intersection between Cayley spectrahedron and the region $x-y\in[-2t,2t]$. Expressing the solutions of P$(^4_t)$ in terms of the variables $x(t),y(t),z(t)$, we have: 
\begin{equation*}
(x(t),y(t),z(t))=
    \begin{cases}
    (0,0,-1)&\text{for $t\in(-1,0]$},
    \\ 
    (t,-t,-1) &\text{for } t\in\left(\frac{1}{2k-1},\frac{1}{2k}\right),\ k=1,2,\dots
    \\ 
    \{(\alpha,-\alpha,-1)\;|\;\alpha\in[-t,t]\} &\text{for $t=\frac{1}{k}$},\qquad \qquad\,\ k=1,2,\dots
    \\
    (-t,t,-1) &\text{for } t\in\left(\frac{1}{2k},\frac{1}{2k+1}\right),\ k=1,2,\dots
    \end{cases}
\end{equation*}  

For every $t\in(-1,0]$, $\mathcal{P}^*(t)$ is continuous and single-valued. The parameter sequence $\{t_k\}_{k=1}^\infty\subseteq(0,1]$ defined by $t_k:=\frac{1}{k} $ is such that $\lim_{k\to\infty}t_k=0$ and at each $t_k$ a loss of inner semi-continuity occurs. Hence, $t=0$ is an irregular accumulation point, according to Def. \ref{def: type_accu}

\begin{figure}[htb!]
    \centering
    \includegraphics[width=0.73\textwidth]{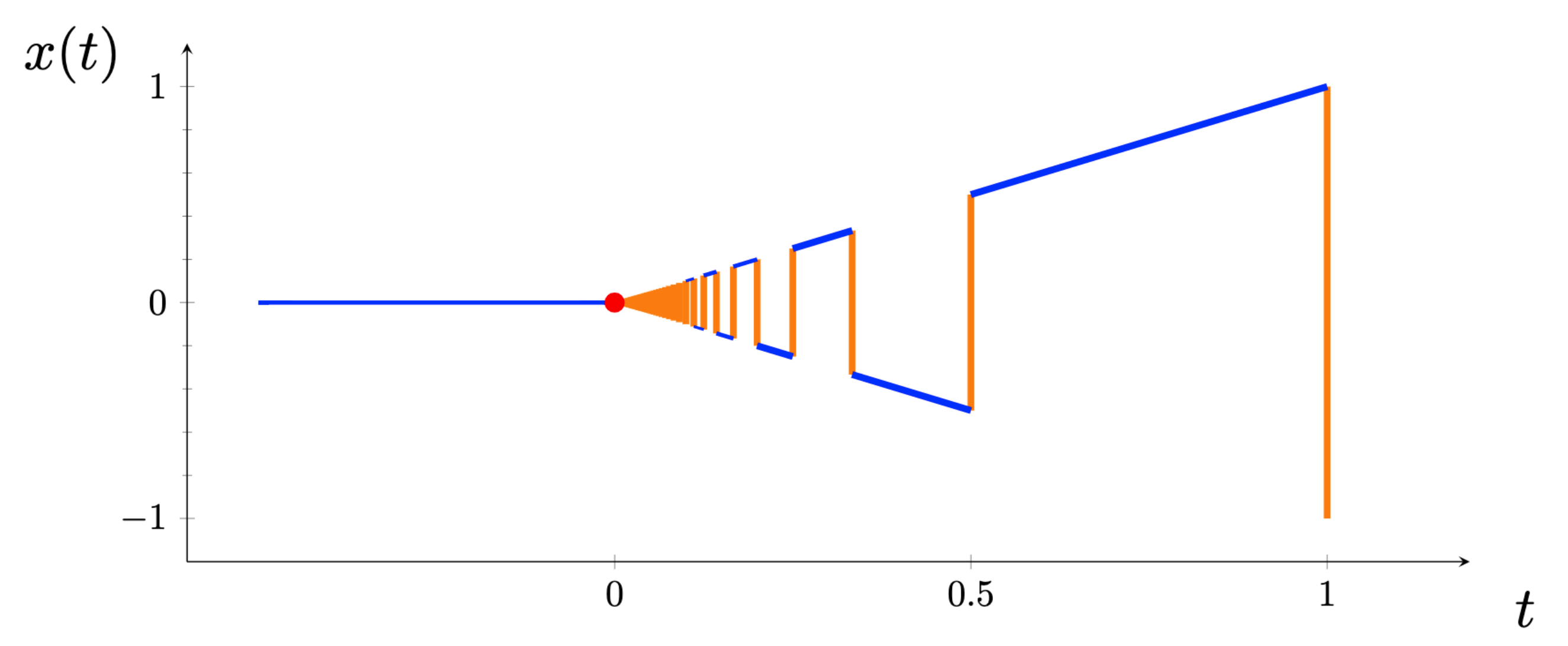}
    \caption{Graph of the $x$ coordinate of the optimal set of $(P^4_t)$ as a function of time $t$. The blue segments correspond to regular points, the red dot corresponds to an irregular accumulation point, and the orange vertical segments correspond to discontinuous isolated multiple-points, where the solution is multiple valued.}
    \label{fig:example_4a}
\end{figure} 

In the following, we also provide an example of an accumulation point for a sequence of continuous bifurcation points.
     
For $t\in(-1,1)$, consider the SDP
\begin{gather}
    \operatorname{min} \ z\nonumber\\
    s.t.\quad\begin{pmatrix}
    1&x&y&0&0\\
    x&1&z&0&0\\
    y&z&1&0&0\\
    0&0&0&2h(t)&x-y\\
    0&0&0&x-y&2h(t)
    \end{pmatrix}\succeq0\tag{P$^5_t$},
\end{gather}
where 
\[
h(t):=\begin{cases}
t\sin^2{\frac{\pi}{t}}&\text{if }t>0,\\
0&\text{otherwise}.
\end{cases}
\]
For $t\le0$ and for $t=1/k,\;k=1,2,\dots$ the feasible region is the intersection between Cayley spectrahedron and the plane $x-y=0$, while for $t\in\left(1/k,1/(k+1)\right),\;k=1,2,\dots$ the feasible region is the intersection between Cayley spectrahedron and the region $x-y\in[-2h(t),2h(t)]$. Writing the solutions of (P$^5_t)$  in terms of the variables $x(t),y(t),z(t)$, we have: 
\begin{equation*}
    (x(t),y(t),z(t))=
    \begin{cases}
        (0,0,-1)&\text{for $t\in(-1,0]$}, 
        \\
        \{(\alpha,-\alpha,-1)\;|\;\alpha\in[-h(t),h(t)]\} &\text{for } t\in\left(\frac{1}{k},\frac{1}{k+1}\right),\ k=1,2,\dots
        \\
        (0,0,-1)&\text{for $t=\frac{1}{k}$},\qquad \quad\ \ k=1,2,\dots
    \end{cases}
\end{equation*} 

For every $t\in(-1,1)$, $\mathcal{P}^*(t)$ is continuous. The parameter sequence $\{t_k\}_{k=1}^\infty\subseteq(0,1]$ defined by $t_k:=\frac{1}{k} $ is such that $\lim_{k\to\infty}t_k=0$ and  each $t_k$ is a continuous bifurcation  point. Hence, $t=0$ is an irregular accumulation point, according to Def. \ref{def: type_accu}.

\begin{figure}[h]
        \centering
        \includegraphics[width=0.73\textwidth]{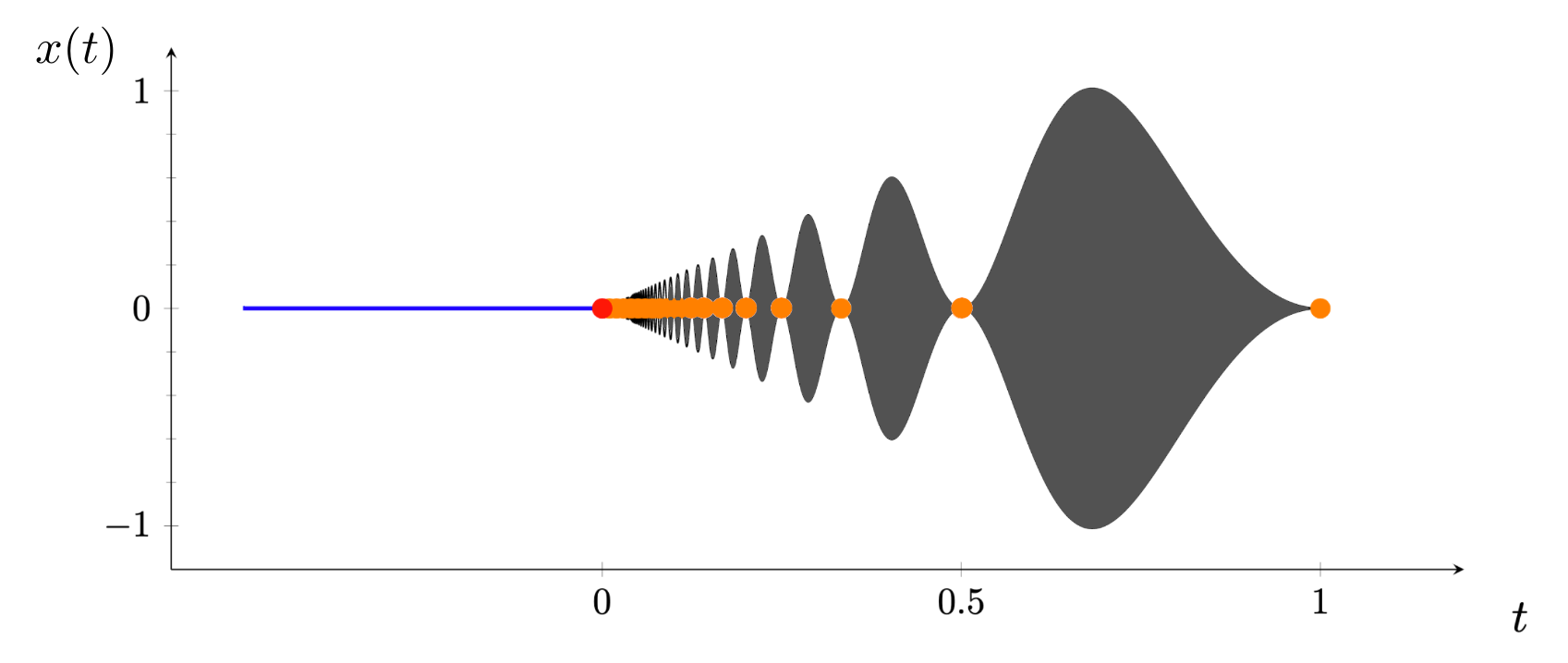}
        \caption{Graph of the $x$ coordinate of the optimal set of (P$^5_t)$ as a function of time $t$. The blue segment consists of regular points, the red dot corresponds to an irregular accumulation point, and the orange dots correspond to continuous bifurcation points. The gray region corresponds to times intervals where the optimal solution is multi-valued.}
        \label{fig:example_4b}
    \end{figure}   
    
\section{Discussion}




Our approach draws upon a long history of work in parametric optimization. In particular, the pioneering work of \cite[Chapter 2]{guddat1990parametric} outlined a classification of solutions to univariate parametric non-linear constrained optimization problems. 
There, precise algebraic conditions are shown for points satisfying first-order optimality conditions to be \textit{non-degenerate} (see Remark \ref{rem: non_degeneracy}). These points exhibit a regular behavior. For \textit{degenerate} points, four different types are defined according to which subset of non-degeneracy conditions is violated. 
Analogously, our approach also starts by considering algebraic conditions that ensure a regular behavior, but our classification of irregular points was made according to the local behavior of the trajectory of solutions at the point considered, rather than according to different sets of algebraic conditions.
\newline

We notice that regular points and discontinuous isolated multiple points, defined as in Definitions \ref{def: type_regular} and \ref{def: type_isol_non_uni} respectively, were first identified by \cite{hauenstein2019computing} (see e.g. Example 1 there) within the optimal partition approach to parametric analysis for linearly  parametrized SDP. Furthermore, non-differentiable points (Definition \ref{def: type_non_diff}) can be easily derived from their results. 
 
Our work can hence be seen as a completion of the effort of \cite{hauenstein2019computing}. 
Likewise, in our analysis, Theorem \ref{thm: finite_bad_result} relies on Theorem \ref{thm: finiteness_of_singularity} and 
Theorem \ref{thm: single_valued_differentiable}.
There, the proof of Theorem \ref{thm: finiteness_of_singularity}  uses the technique of \cite{hauenstein2019computing}, while Theorem \ref{thm: single_valued_differentiable} is essentially an application of the implicit function theorem, implying that this can be applied almost everywhere. Theorem \ref{thm: main_result} suggests that when, instead, the assumptions for implicit function theorem do not hold almost everywhere, this allows for a broader range of possible behaviors, listed in the last row of Table \ref{tab: tvsdp_behaviors}. 
\newline

\begin{table}[htb!]
    \centering
    \footnotesize
    \begin{tabular}{l|l}
        \hline
        \textbf{Problem assumptions}&\textbf{Type of points}
        \\
        \hline
        SDP with LICQ, continuous data, &Regular points\\
        \phantom{SDP with }strict feasibility, and a \textit{non-singular} time & Non-differentiable points\\
        &Discontinuous isolated multiple points
        \\
        \hline
        SDP with LICQ,  continuous data,  & Regular points\\
        \phantom{SDP with }strict feasibility, without a \textit{non-singular} time & Non-differentiable points
        \\& Discontinuous isolated multiple points
        \\&Discontinuous non-isolated multiple points\\
        & Continuous bifurcation points\\
        & Irregular accumulation points
        \\ 
        \hline
    \end{tabular} 
    \vspace{.3cm}
    \caption{Assumptions on {parametric} SDP and associated possible type of points}
    \label{tab: tvsdp_behaviors}
\end{table}
From the point of view of formulating a {parametric} SDP, the key insight of \cite{hauenstein2019computing} and ours is that even seemingly strong and standard assumptions such as the continuity of the data and primal-dual strict feasibility are not sufficient to prevent pathological behavior. We presented a complete characterization of such behaviors. Thereby, we showed that guaranteeing the existence of a generic non-singular point along the trajectory suffices to prevent highly pathological behaviors. However, this does not prevent from a finite number of losses of differentiability or isolated losses of uniqueness to occur.
\newline 

One may also be interested in understanding how the main result of this paper specializes to restricted classes of {parametric} SDP, such as {parametric}  linear programming (LP) and {parametric}  second order cone programming (SOCP). In the first case, if the data are assumed to be continuous functions, one can easily construct an example of each type of behaviors of the trajectory of solutions described in Definitions 12-17.
\newline 

For example, for $t\in(-1,1)$ consider:

\begin{enumerate}
    \item $\min\{x :x\ge 1+t\}$.
    \item $\min\{x :x\ge|t|\}$.
    \item $\min\{tx :-1\le x\le1\}$.
    \item $\min\{f(t)x :-1\le x\le1\}$, with $ f(t)=t$ if $ t\le0$, otherwise $f(t)=0 $.
    \item $\min\{0x :-g(t)\le x\le g(t)\}$, with $ g(t)=0$ if $ t\le0$, otherwise $g(t)=t $.
    \item $\min\{0x :-h(t)\le x\le h(t)\}$, with $h(t)=0$ if $t>0$, otherwise $h(t)=t\sin^2{\frac{\pi}{t}}$.
\end{enumerate}  
At $t^*=0$, $x^*=0$ is 1. a regular point, 2. a non-differentiable point, 3. an isolated discontinuous multiple point, 4. a non-isolated discontinuous multiple point, 5 a continuous bifurcation point, 6 an irregular accumulation point.
Hence, restricting to the class of {parametric} LP does not exclude any type of point. It follows that also in the case of SOCP, a class that generalize LP, all the type of points can possibly appear. From this point of view, it is surprising that the trajectories of solution to {parametric} SDP, a class of optimization problems much wider than LP, does not present, in the general framework that we adopted, any behavior which does not already show up in {parametric} LP. However, we believe that under a set of assumptions more specific than the one that we adopted in Theorem \ref{thm: main_result}, some type of behaviors may be ruled out in {parametric} LP, but not in {parametric} SDP. 
Take as an example non-differentiable points (see Def. \ref{def: type_non_diff}). If one assumes that the time dependence of the data is smooth, e.g. polynomial, non differentiable points can still appear in {parametric} SDP (see the first example of Section 4). This is due to the facial geometry of SDP, where positively curved surfaces appear, which must then entirely consist of extreme points (0-dimensional faces). Instead, in LP, extreme points are always isolated, so that when the solution is unique, this must coincide with a fixed extreme point. If the time dependence is smooth, this should imply that the feasible set, hence its extreme points, should also move smoothly, preventing non-differentiable points to show up.  
The investigation of such distinctions between {parametric} LP and {parametric} SDP may be an interesting direction for future research.

\section{Conclusion}

We used set-valued analysis to describe and study the trajectory of solutions to {parametric} SDP. 
The analysis we carried out brought us to define six different types of points, according to the local structure of the solutions trajectory.
Our main result consists in proving that under standard assumptions, there are no other types of points.
\newline 

One could extend our research by weakening our assumptions:  continuity of the data dependence on the parameter, and primal and dual strict feasibility throughout the parameterization interval. These requirements avoid highly degenerate situations. In particular, without continuity of the data, one can expect the trajectory to potentially present a lot of irregularities, e.g., it may fail to be both inner and outer semi-continuous, while, as Theorem \ref{thm: optimal_sets_osc} shows, under the continuity of the data outer semi-continuity is ensured.
When strict feasibility is lost, two additional forms of degenerate behavior might occur: the optimal value may not be attained at any feasible point, or there may be a strictly positive duality gap between the primal and dual optimal values.
It is not clear whether there could be other types too, perhaps akin to irregular accumulation points. 
\newline 
 
Finally, the properties of specific classes of trajectories of solutions in specific applications may be of considerable interest. 

\section*{Acknowledgments.}

This research has been supported by the OP RDE funded project 
CZ.02.1.01/0.0/0.0/16\_019/0000765 ``Research Center for Informatics''.
\newpage
\phantom{\cite{drusvyatskiy2017many}} 
 
\bibliography{bib} 

\end{document}